\documentclass[twoside,11pt]{article}
\title{} \author{} \date{}
\usepackage{latexsym,amssymb,times}
\input amssym.def
\newtheorem{te}{Theorem}[section]

\newtheorem{fac}[te]{Fact}

\newtheorem{cla}[te]{Claim}
\newtheorem{rem}[te]{Remark}
\newtheorem{ex}[te]{Example}

\begin{document}
\thispagestyle{plain}
\begin{center}
           {\large \bf {MAXIMALLY EMBEDDABLE COMPONENTS}}
\end{center}
\begin{center}
{\small\bf Milo\v s S.\ Kurili\'c}\\
         {\small Department of Mathematics and Informatics, University of Novi Sad, \\
         Trg Dositeja Obradovi\'ca 4, 21000 Novi Sad, Serbia.\\
                                     e-mail: milos@dmi.uns.ac.rs}
\end{center}
\begin{abstract}
\noindent
We investigate the partial orderings of the form $\langle  {\mathbb P} ({\mathbb X} ), \subset\rangle $, where ${\mathbb X} =\langle  X,\rho \rangle $ is a countable binary relational structure and
${\mathbb P} ({\mathbb X} )$ the set of the domains of its isomorphic substructures and show that if the components of ${\mathbb X}$ are maximally embeddable and satisfy an
additional condition related to connectivity, then the poset $\langle  {\mathbb P} ({\mathbb X} ), \subset\rangle $ is forcing equivalent to
a finite power of $(P(\omega )/\mathop{\rm Fin}\nolimits  )^+$,  or to $(P(\omega \times \omega)/(\mathop{\rm Fin}\nolimits  \times \mathop{\rm Fin}\nolimits  ))^+$, or to the direct product
$(P(\Delta )/{\mathcal E}{\mathcal D}_{\mathrm{fin}}  )^+  \times  ((P(\omega )/\mathop{\rm Fin}\nolimits )^+)^n$, for some $n\in\omega$. In particular we obtain forcing equivalents
of the posets of copies of countable equivalence relations, disconnected ultrahomogeneous graphs and some partial orderings.


\end{abstract}
\section{Introduction}\label{S4000}
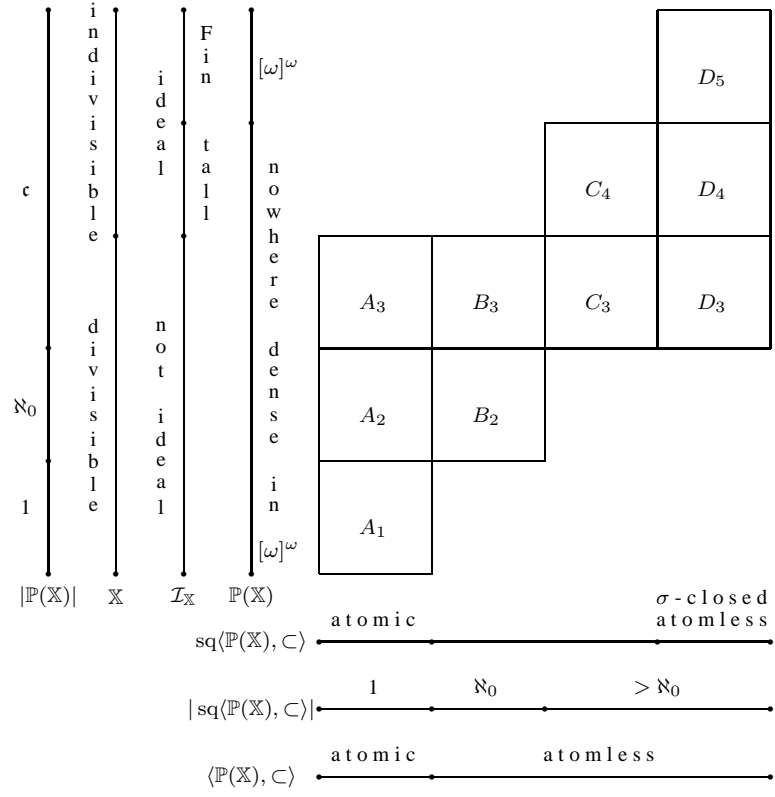
\begin{figure}[htb]
\begin{center}
\unitlength 0.6mm 
\linethickness{0.4pt}
\ifx\plotpoint\undefined\newsavebox{\plotpoint}\fi 


\begin{picture}(180,190)(0,0)


\put(70,10){\line(1,0){100}}
\put(70,25){\line(1,0){100}}
\put(70,40){\line(1,0){100}}
\put(10,55){\line(0,1){125}}
\put(25,55){\line(0,1){125}}
\put(40,55){\line(0,1){125}}
\put(55,55){\line(0,1){125}}
\put(70,55){\line(0,1){75}}
\put(95,55){\line(0,1){75}}
\put(120,80){\line(0,1){50}}
\put(120,130){\line(0,1){25}}
\put(145,105){\line(0,1){75}}
\put(170,105){\line(0,1){75}}
\put(70,55){\line(1,0){25}}
\put(70,80){\line(1,0){50}}
\put(70,105){\line(1,0){100}}
\put(70,130){\line(1,0){100}}
\put(120,155){\line(1,0){50}}
\put(145,180){\line(1,0){25}}


\put(70,10){\circle*{1}}
\put(95,10){\circle*{1}}
\put(170,10){\circle*{1}}
\put(70,25){\circle*{1}}
\put(95,25){\circle*{1}}
\put(120,25){\circle*{1}}
\put(170,25){\circle*{1}}
\put(70,40){\circle*{1}}
\put(145,40){\circle*{1}}
\put(170,40){\circle*{1}}
\put(10,55){\circle*{1}}
\put(10,80){\circle*{1}}
\put(10,105){\circle*{1}}
\put(10,180){\circle*{1}}
\put(25,55){\circle*{1}}
\put(25,130){\circle*{1}}
\put(25,180){\circle*{1}}
\put(40,55){\circle*{1}}
\put(40,130){\circle*{1}}
\put(40,155){\circle*{1}}
\put(40,180){\circle*{1}}
\put(55,55){\circle*{1}}
\put(55,155){\circle*{1}}
\put(55,180){\circle*{1}}
\put(95,40){\circle*{1}}

\scriptsize
\put(82,15){\makebox(0,0)[cc]{a t o m i c}}%
\put(132,15){\makebox(0,0)[cc]{a t o m l e s s}}%
\put(82,45){\makebox(0,0)[cc]{a t o m i c}}%

\put(82,30){\makebox(0,0)[cc]{1}}%
\put(107,30){\makebox(0,0)[cc]{$\aleph _0$}}%
\put(145,30){\makebox(0,0)[cc]{$>\aleph _0$}}%

\put(157,50){\makebox(0,0)[cc]{$\sigma$ - c l o s e d}}%
\put(157,45){\makebox(0,0)[cc]{a t o m l e s s}}%

\put(5,70){\makebox(0,0)[cc]{1}}%
\put(5,92){\makebox(0,0)[cc]{$\aleph _0$}}%
\put(5,140){\makebox(0,0)[cc]{${\mathfrak c}$}}%
\put(20,180){\makebox(0,0)[cc]{i}}%
\put(20,175){\makebox(0,0)[cc]{n}}%
\put(20,170){\makebox(0,0)[cc]{d}}%
\put(20,165){\makebox(0,0)[cc]{i}}%
\put(20,160){\makebox(0,0)[cc]{v}}%
\put(20,155){\makebox(0,0)[cc]{i}}%
\put(20,150){\makebox(0,0)[cc]{s}}%
\put(20,145){\makebox(0,0)[cc]{i}}%
\put(20,140){\makebox(0,0)[cc]{b}}%
\put(20,135){\makebox(0,0)[cc]{l}}%
\put(20,130){\makebox(0,0)[cc]{e}}%
\put(20,110){\makebox(0,0)[cc]{d}}%
\put(20,105){\makebox(0,0)[cc]{i}}%
\put(20,100){\makebox(0,0)[cc]{v}}%
\put(20,95){\makebox(0,0)[cc]{i}}%
\put(20,90){\makebox(0,0)[cc]{s}}%
\put(20,85){\makebox(0,0)[cc]{i}}%
\put(20,80){\makebox(0,0)[cc]{b}}%
\put(20,75){\makebox(0,0)[cc]{l}}%
\put(20,70){\makebox(0,0)[cc]{e}}%
\put(35,165){\makebox(0,0)[cc]{i}}%
\put(35,160){\makebox(0,0)[cc]{d}}%
\put(35,155){\makebox(0,0)[cc]{e}}%
\put(35,150){\makebox(0,0)[cc]{a}}%
\put(35,145){\makebox(0,0)[cc]{l}}%
\put(35,110){\makebox(0,0)[cc]{n}}%
\put(35,105){\makebox(0,0)[cc]{o}}%
\put(35,100){\makebox(0,0)[cc]{t}}%
\put(35,90){\makebox(0,0)[cc]{i}}%
\put(35,85){\makebox(0,0)[cc]{d}}%
\put(35,80){\makebox(0,0)[cc]{e}}%
\put(35,75){\makebox(0,0)[cc]{a}}%
\put(35,70){\makebox(0,0)[cc]{l}}%
\put(45,175){\makebox(0,0)[cc]{F}}%
\put(45,170){\makebox(0,0)[cc]{i}}%
\put(45,165){\makebox(0,0)[cc]{n}}%
\put(45,150){\makebox(0,0)[cc]{t}}%
\put(45,145){\makebox(0,0)[cc]{a}}%
\put(45,140){\makebox(0,0)[cc]{l}}%
\put(45,135){\makebox(0,0)[cc]{l}}%
\put(61,167){\makebox(0,0)[cc]{$[\omega ]^{\omega }$}}%
\put(60,145){\makebox(0,0)[cc]{n}}%
\put(60,140){\makebox(0,0)[cc]{o}}%
\put(60,135){\makebox(0,0)[cc]{w}}%
\put(60,130){\makebox(0,0)[cc]{h}}%
\put(60,125){\makebox(0,0)[cc]{e}}%
\put(60,120){\makebox(0,0)[cc]{r}}%
\put(60,115){\makebox(0,0)[cc]{e}}%
\put(60,105){\makebox(0,0)[cc]{d}}%
\put(60,100){\makebox(0,0)[cc]{e}}%
\put(60,95){\makebox(0,0)[cc]{n}}%
\put(60,90){\makebox(0,0)[cc]{s}}%
\put(60,85){\makebox(0,0)[cc]{e}}%
\put(60,75){\makebox(0,0)[cc]{i}}%
\put(60,70){\makebox(0,0)[cc]{n}}%
\put(61,60){\makebox(0,0)[cc]{$[\omega ]^{\omega }$}}%



\put(10,50){\makebox(0,0)[cc]{$|{\mathbb P} ({\mathbb X} )|$}}%
\put(25,50){\makebox(0,0)[cc]{${\mathbb X}$}}%
\put(40,50){\makebox(0,0)[cc]{${\mathcal I}_{\mathbb X}$}}%
\put(55,50){\makebox(0,0)[cc]{${\mathbb P}({\mathbb X})$}}%
\put(55,40){\makebox(0,0)[cc]{$\mathop{\rm sq}\nolimits  \langle {\mathbb P} ({\mathbb X} ) , \subset \rangle$}}%
\put(55,25){\makebox(0,0)[cc]{$|\mathop{\rm sq}\nolimits  \langle {\mathbb P} ({\mathbb X} ) , \subset \rangle|$}}%
\put(55,10){\makebox(0,0)[cc]{$\langle {\mathbb P} ({\mathbb X} ) , \subset \rangle$}}%


\put(82,65){\makebox(0,0)[cc]{$A_1$}}%
\put(82,90){\makebox(0,0)[cc]{$A_2$}}%
\put(82,115){\makebox(0,0)[cc]{$A_3$}}%
\put(107,90){\makebox(0,0)[cc]{$B_2$}}%
\put(107,115){\makebox(0,0)[cc]{$B_3$}}%
\put(132,115){\makebox(0,0)[cc]{$C_3$}}%
\put(132,140){\makebox(0,0)[cc]{$C_4$}}%
\put(157,115){\makebox(0,0)[cc]{$D_3$}}%
\put(157,140){\makebox(0,0)[cc]{$D_4$}}%
\put(157,165){\makebox(0,0)[cc]{$D_5$}}%
\end{picture}

\end{center}

\vspace{-7mm}

\caption{Binary relations on countable sets}\label{F4001}
\end{figure}

The posets of the form $\langle  {\mathbb P} ({\mathbb X} ), \subset\rangle $, where ${\mathbb X} $ is a relational structure and
${\mathbb P} ({\mathbb X} )$ the set of the domains of its isomorphic substructures, were investigated in \cite{Ktow}.
In particular, a classification of countable binary structures
related to the forcing-related properties of the posets of their copies is described in Diagram \ref{F4001}:
for the structures from column $A$ (resp.\ $B$; $D$) the corresponding posets are forcing equivalent to the trivial poset
(resp.\ the Cohen forcing, $\langle  {}^{<\omega }2, \supset\rangle $;
an $\omega _1 $-closed atomless poset) and, for the structures from the class $C_4$, the posets of copies are forcing equivalent to the quotients
of the form $P(\omega )/{\mathcal I}$, for some co-analytic tall ideal ${\mathcal I}$.

The aim of the paper is to
investigate a subclass of column $D$, the class of structures ${\mathbb X}$ for which  the separative quotient $\mathop{\rm sq}\nolimits  \langle  {\mathbb P} ({\mathbb X} ), \subset\rangle $ is
an $\omega _1$-closed and atomless poset (containing, for example, the class of all countable scattered linear orders \cite{Kurscatt}).
Clearly, such a classification
depends on the model of
set theory in which we work. For example, under the CH all the structures from column $D$ are in the same class
(having the posets of copies forcing equivalent to the algebra $P(\omega )/\!\mathop{\rm Fin}\nolimits $ without zero), but this is not true
in, for example, the Mathias model.

Applying the main theorem of the paper, proved in Section \ref{S4006}, in Section \ref{S4008} we obtain forcing equivalents
of the posets of copies of countable equivalence relations, disconnected ultrahomogeneous graphs and some partial orderings.

\section{Preliminaries}\label{S4001}
Let ${\mathbb P} =\langle  P , \leq \rangle $ be a pre-order. Then
$p\in P$ is an {\it atom}, in notation  $p\in \mathop{\rm At}\nolimits ({\mathbb P} )$, iff each  $q,r\leq p$ are compatible (there is $s\leq q,r$).
${\mathbb P} $ is called {\it atomless} iff $\mathop{\rm At}\nolimits ({\mathbb P} )=\emptyset$; {\it atomic} iff $\mathop{\rm At}\nolimits ({\mathbb P} )$ is dense in ${\mathbb P}$.
If $\kappa $ is a regular cardinal, ${\mathbb P} $ is called $\kappa ${\it -closed} iff for each
$\gamma <\kappa $ each sequence $\langle  p_\alpha :\alpha <\gamma\rangle $ in $P$, such that $\alpha <\beta \Rightarrow p_{\beta}\leq p_\alpha $,
has a lower bound in $P$. $\omega _1$-closed pre-orders are called {\it $\sigma$-closed}.
Two pre-orders ${\mathbb P}$ and ${\mathbb Q}$ are called {\it forcing equivalent} iff they produce the same generic extensions.

A partial order ${\mathbb P} =\langle  P , \leq \rangle $ is called
{\it separative} iff for each $p,q\in P$ satisfying $p\not\leq q$ there is $r\leq p$ such that $r \perp q$.
The {\it separative modification} of ${\mathbb P}$
is the separative pre-order $\mathop{\rm sm}\nolimits ({\mathbb P} )=\langle  P , \leq ^*\rangle $, where
$p\leq ^* q \Leftrightarrow \forall r\leq p \; \exists s \leq r \; s\leq q $.
The {\it separative quotient} of ${\mathbb P}$
is the separative partial order $\mathop{\rm sq}\nolimits  ({\mathbb P} )=\langle P /\!\! =^* , \trianglelefteq \rangle$, where
$p = ^* q \Leftrightarrow p \leq ^* q \land q \leq ^* p\;$ and $\;[p] \trianglelefteq [q] \Leftrightarrow p \leq ^* q $.

Let $\mathop{\rm Fin}\nolimits  = [\omega ]^{<\omega}$ and $\Delta =\{ \langle  m,n\rangle \in {\mathbb N} \times {\mathbb N} : n\leq m\}$. Then the ideals $\mathop{\rm Fin}\nolimits  \times \mathop{\rm Fin}\nolimits  \subset P(\omega\times \omega )$
and ${\mathcal E}{\mathcal D}_{\mathrm{fin}} \subset P(\Delta)$ are defined by:

$\mathop{\rm Fin}\nolimits  \times \mathop{\rm Fin}\nolimits  = \{ S\subset \omega \times \omega : \exists j\in \omega \; \forall i\geq j \;|S\cap (\{ i\}\times \omega )|<\omega  \}$ and

${\mathcal E}{\mathcal D}_{\mathrm{fin}}  = \{ S\subset \Delta : \exists r\in {\mathbb N} \;\; \forall m \in {\mathbb N} \;\; |S\cap (\{ m \} \times \{ 1,2,\dots ,m\})|\leq r  \}$.

\noindent
By ${\mathfrak h}({\mathbb P} )$ we denote the {\it distributivity number} of a poset ${\mathbb P}$. In particular, for $n\in {\mathbb N}$, let
${\mathfrak h}_n={\mathfrak h}(((P(\omega )/\mathop{\rm Fin}\nolimits  )^+ )^n )$; thus ${\mathfrak h}={\mathfrak h}_1$. The following statements will be used in the paper.
\begin{fac}    \rm \label{T4043}
(Folklore) If ${\mathbb P} _i$, $i\in I$, are $\kappa $-closed pre-orders, then $\prod _{i\in I}{\mathbb P} _i$ is $\kappa $-closed.
\end{fac}
\begin{fac}    \rm \label{T4042}
(Folklore) Let ${\mathbb P} , {\mathbb Q} $ and ${\mathbb P} _i$, $i\in I$, be partial orderings. Then

(a) ${\mathbb P}$, $\mathop{\rm sm}\nolimits ({\mathbb P})$ and $\mathop{\rm sq}\nolimits  ({\mathbb P})$ are forcing equivalent forcing notions;

(b) ${\mathbb P}$ is atomless iff $\mathop{\rm sm}\nolimits ({\mathbb P} )$ is atomless iff $\mathop{\rm sq}\nolimits  ({\mathbb P} )$ is atomless;

(c) $\mathop{\rm sm}\nolimits ({\mathbb P} )$ is $\kappa $-closed iff $\mathop{\rm sq}\nolimits  ({\mathbb P} )$ is $\kappa $-closed;

(d) ${\mathbb P} \cong {\mathbb Q}$ implies that $\mathop{\rm sm}\nolimits {\mathbb P} \cong  \mathop{\rm sm}\nolimits {\mathbb Q}$ and $\mathop{\rm sq}\nolimits  {\mathbb P} \cong  \mathop{\rm sq}\nolimits  {\mathbb Q}$;

(e) $\mathop{\rm sm}\nolimits (\prod _{i\in I}{\mathbb P} _i) = \prod _{i\in I}\mathop{\rm sm}\nolimits {\mathbb P} _i$;

(f) $\mathop{\rm sq}\nolimits  (\prod _{i\in I}{\mathbb P} _i) \cong \prod _{i\in I}\mathop{\rm sq}\nolimits  {\mathbb P} _i$.
\end{fac}
\begin{fac}   \rm \label{T4056}
(Folklore) Let ${\mathbb P}$ be an atomless separative pre-order. Then we have

(a) If $\omega _1={\mathfrak c}$ and ${\mathbb P}$ is $\omega _1$-closed of size ${\mathfrak c}$, then ${\mathbb P}$ is forcing equivalent to $(\mathop{\rm Coll}\nolimits (\omega _1 , \omega _1 ))^+$
or, equivalently, to $(P(\omega )/\mathop{\rm Fin}\nolimits )^+$;

(b) If ${\mathfrak t}={\mathfrak c}$ and ${\mathbb P}$ is ${\mathfrak t}$-closed of size ${\mathfrak t}$, then ${\mathbb P}$ is forcing equivalent to
$(\mathop{\rm Coll}\nolimits ({\mathfrak t},{\mathfrak t}))^+$
or, equivalently, to $(P(\omega )/\mathop{\rm Fin}\nolimits )^+$.
\end{fac}
\begin{fac}    \rm  \label{T4091}
(a) $\mathop{\rm sm}\nolimits (\langle  [\omega ]^\omega , \subset \rangle ^n ) =\langle  [\omega ]^\omega , \subset ^*\rangle ^n$ and
$\mathop{\rm sq}\nolimits  (\langle  [\omega ]^\omega , \subset \rangle ^n) =((P(\omega )/\mathop{\rm Fin}\nolimits  )^+ )^n$ are forcing equivalent,
${\mathfrak t}$-closed atomless  pre-orders of size ${\mathfrak c}$.

(b) {\bf (Shelah and Spinas \cite{SheSpi1})} Con(${\mathfrak h}_{n+1}  < {\mathfrak h}_n$), for each $n\in {\mathbb N}$.

(c) {\bf (Szyma\'nski and Zhou \cite{Szym})} $(P(\omega \times \omega)/(\mathop{\rm Fin}\nolimits  \times \mathop{\rm Fin}\nolimits  ))^+$
is an $\omega _1$-closed, but not $\omega _2$-closed atomless  poset.

(d) {\bf (Hern\'andez-Hern\'andez \cite{Her})}
Con(${\mathfrak h}((P(\omega \times \omega )/(\mathop{\rm Fin}\nolimits  \times \mathop{\rm Fin}\nolimits  ))^+) < {\mathfrak h}$).

(e) {\bf (Brendle \cite{Bren})}
Con(${\mathfrak h}(( P(\Delta )/{\mathcal E}{\mathcal D}_{\mathrm{fin}}  )^+) <{\mathfrak h}$).
\end{fac}
\begin{fac}    \rm  \label{T4090}
If $\langle  P , \leq_P \rangle $ and $\langle  Q , \leq_Q \rangle $ are  partial orderings and $f: P\rightarrow Q$, where

(i) $\forall p_1 , p_2 \in P \; (p_1 \leq _P p_2 \Rightarrow f(p_1) \leq _Q f(p_2))$,

(ii) $\forall p_1 , p_2 \in P \; (p_1 \perp _P p_2 \Rightarrow f(p_1) \perp _Q f(p_2))$,

(iii) $f[P]=Q$,

\noindent
then $\mathop{\rm sq}\nolimits  {\mathbb P} \cong \mathop{\rm sq}\nolimits  {\mathbb Q}$.
\end{fac}
\noindent{\bf Proof.}
We have $\mathop{\rm sm}\nolimits {\mathbb P} = \langle  P , \leq_P^* \rangle $, $\mathop{\rm sq}\nolimits  {\mathbb P} = \langle  P/\!\!=_P , \trianglelefteq _P \rangle $,
$\mathop{\rm sm}\nolimits {\mathbb Q} = \langle  Q , \leq_Q^* \rangle $ and $\mathop{\rm sq}\nolimits  {\mathbb Q} = \langle  Q/\!\!=_Q , \trianglelefteq _Q \rangle $,
where for each $p_1,p_2\in P$ and each $q_1,q_2\in Q$
\begin{equation}\label{EQ4090}
p_1 \leq_P^* p_2 \Leftrightarrow \forall p\leq _P p_1 \; \exists p'\leq _P p,p_2 ,
\end{equation}
\begin{equation}\label{EQ4091}
p_1 =_P p_2 \Leftrightarrow p_1 \leq_P^* p_2 \land p_2 \leq_P^* p_1 \;\; \mbox{ and }\;\;
[p_1] {\trianglelefteq _P} [p_2]\Leftrightarrow p_1 \leq_P^* p_2 ,
\end{equation}
\begin{equation}\label{EQ4094}
q_1 \leq_Q^* q_2 \Leftrightarrow \forall q\leq _Q q_1 \; \exists q'\leq _Q q,q_2 ,
\end{equation}
\begin{equation}\label{EQ4095}
q_1 =_Q q_2 \Leftrightarrow q_1 \leq_Q^* q_2 \land q_2 \leq_Q^* q_1 \;\; \mbox{ and }\;\;
[q_1] {\trianglelefteq _Q} [q_2]\Leftrightarrow q_1 \leq_Q^* q_2.
\end{equation}
{\it Claim.} $p_1 \leq _P^* p_2 \Leftrightarrow f(p_1) \leq _Q^* f(p_2)$, for each $p_1,p_2\in P$.

\vspace{2mm}
\noindent
{\it Proof of Claim.}
($\Rightarrow$) Let $p_1 \leq _P^* p_2$. According to (\ref{EQ4094}) we prove
\begin{equation}\label{EQ4093}
\forall q\leq _Q f(p_1) \; \exists q'\leq _Q q,f(p_2) .
\end{equation}
If $q\leq _Q f(p_1)$ then, by (iii) there is $p_3\in P$ such that $f(p_3)=q$. By (ii) and since $f(p_3)\leq _Q f(p_1)$, there is
$p_4 \leq _P p_3 , p_1$ and, by (\ref{EQ4090}), there is $p_5\leq _P p_4, p_2$, which, by (i), implies $f(p_5)\leq _Q f(p_2)$. Since $p_5 \leq _P p_4 \leq _P p_3$
by (i) we have $f(p_5)\leq _Q f(p_3)=q$ and $q'=f(p_5)$ satisfies (\ref{EQ4093}).

($\Leftarrow$) Assuming (\ref{EQ4093}) we prove that $p_1 \leq _P^* p_2$. If $p\leq _P p_1$, then, by (i), $f(p)\leq _Q f(p_1)$
and, by (\ref{EQ4093}), there is $q'\leq _Q f(p),f(p_2)$ and, by (ii), there is $p'\leq _P p,p_2$ and Claim is proved.

\vspace{2mm}
\noindent
Now we show that $\langle  P/\!\!=_P , \trianglelefteq _P \rangle \cong _F \langle  Q/\!\!=_Q , \trianglelefteq _Q \rangle $, where
$F([p])=[f(p)]$.

By Claim, (\ref{EQ4091}) and (\ref{EQ4095}), for each $p_1,p_2\in P$ we have
$[p_1]=[p_2]$
iff $p_1 =_P p_2$
iff $p_1 \leq_P^* p_2 \land p_2 \leq_P^* p_1$
iff $f(p_1) \leq_Q^* f(p_2) \land f(p_2) \leq_Q^* f(p_1)$
iff $f(p_1) =_Q f(p_2)$
iff $[f(p_1)]=[f(p_2)]$
iff $F([p_1])= F([p_2])$
and $F$ is a well defined injection.
By (iii), for $[q]\in Q/\!\!=_Q$ there is $p\in P$ such that $q=f(p) $. Thus $F([p])=[f(p)]= [q]$ and $F$ is a surjection.

By Claim, (\ref{EQ4091}) and (\ref{EQ4095}) again,
$[p_1] \trianglelefteq _P [p_2]$
iff $p_1 \leq_P^* p_2$
iff $f(p_1) \leq_Q^* f(p_2) $
iff $[f(p_1)] \trianglelefteq _Q [f(p_2)]$
iff $F([p_1])\trianglelefteq _Q F([p_2])$.
Thus $F$ is an  isomorphism.
\hfill $\Box$
\section{Structures and posets of their copies}\label{S4002}
Let $L=\{ R \}$ be a relational language, where $\mathop{\rm ar}\nolimits (R)=2$.
An $L$-structure ${\mathbb X} = \langle  X,  \rho \rangle $ is called a {\it countable structure} iff $|X|=\omega$.
If  $A\subset X$, then $\langle  A, \rho _A \rangle $ is a {\it substructure}
of ${\mathbb X}$, where $\rho _A = \rho  \cap A^2$.
If ${\mathbb Y} =\langle  Y,  \tau   \rangle $ is an $L$-structure too, a map
$f:X \rightarrow Y$ is called an {\it embedding} (we write
${\mathbb X} \hookrightarrow _f {\mathbb Y}$) iff it is an injection and
$\langle  x_1, x_2\rangle \in \rho \Leftrightarrow
\langle  f(x_1), f(x_2)\rangle \in \tau $, for each $\langle  x_1, x_2\rangle \in X^2$.
If ${\mathbb X}$ embeds in ${\mathbb Y}$ we write ${\mathbb X} \hookrightarrow {\mathbb Y}$. Let
$\mathop{\rm Emb}\nolimits ({\mathbb X} , {\mathbb Y} )  =  \{ f: {\mathbb X} \hookrightarrow _f {\mathbb Y} \} $ and, in particular,
$\mathop{\rm Emb}\nolimits ({\mathbb X} )  = \{ f: {\mathbb X} \hookrightarrow _f {\mathbb X} \}$.
If, in addition, $f$ is a surjection, it is an {\it isomorphism} (we write
${\mathbb X} \cong _f {\mathbb Y}$) and the structures ${\mathbb X}$ and ${\mathbb Y}$ are {\it isomorphic},
in notation ${\mathbb X} \cong {\mathbb Y}$.
${\mathbb X}$ and ${\mathbb Y}$ are {\it equimorphic} iff ${\mathbb X}\hookrightarrow {\mathbb Y}$ and ${\mathbb Y}\hookrightarrow {\mathbb X}$.
According to \cite{Fra} a relational structure ${\mathbb X} $ is: {\it $p$-monomorphic} iff
all its substructures of size $p$ are isomorphic; {\it indivisible} iff for each partition
$X=A \cup B$ we have ${\mathbb X} \hookrightarrow A$ or ${\mathbb X}\hookrightarrow B$.

If ${\mathbb X} _i=\langle  X_i, \rho _i \rangle $, $i\in I$, are $L$-structures and $X_i \cap X_j =\emptyset$, for $i\neq j$, then the structure $\bigcup _{i\in I} {\mathbb X} _i =\langle  \bigcup _{i\in I} X_i , \bigcup _{i\in I} \rho _i\rangle $ is the {\it union} of the structures ${\mathbb X} _i$, $i\in I$.

Let $\langle  X,\rho \rangle $ be an $L$-structure and $\rho _{rst}$ the minimal equivalence relation on $X$ containing $\rho $
(the transitive closure  of the relation $\rho _{rs}=\Delta _X \cup \rho \cup \rho ^{-1} $ given by $x \;\rho _{rst} \;y$ iff there are $n\in {\mathbb N}$ and
$z_0 =x , z_1, \dots ,z_n =y$ such that $z_i \;\rho _{rs} \;z_{i+1}$, for each $i<n$).
For $x\in X$ the corresponding equivalence class will be denoted by $[x]$  and called the
{\it component} of $\langle  X,\rho \rangle $ containing $x$. The structure $\langle  X,\rho \rangle $ will be called {\it connected}
iff it has only one component. It is easy to prove (see \cite{Ktow}) that
$\langle  X,\rho \rangle= \langle  \bigcup _{x\in X}[x],\bigcup _{x\in X}\rho _{[x]} \rangle $
is the unique representation of $\langle  X,\rho \rangle $ as a disjoint union of connected relations.

Here we investigate the partial orders
of the form $\langle  {\mathbb P} ({\mathbb X} ), \subset \rangle $, where ${\mathbb X} =\langle  X, \rho \rangle $ is an $L$-structure and
${\mathbb P} ({\mathbb X} )$ the set of its isomorphic substructures, that is
$$
{\mathbb P} ({\mathbb X} )  =  \{ A\subset X : \langle  A, \rho _A  \rangle \cong {\mathbb X}\}
          =  \{ f[X] : f \in \mathop{\rm Emb}\nolimits ({\mathbb X} )\} .
$$
More generally, if  ${\mathbb X} = \langle  X, \rho \rangle $ and ${\mathbb Y} =\langle  Y,  \tau  \rangle $ are two $L$-structures we define
${\mathbb P} ({\mathbb X} , {\mathbb Y}  )  =  \{ B\subset Y : \langle  B,  \tau  _B  \rangle \cong {\mathbb X} \} =  \{ f[X] : f \in \mathop{\rm Emb}\nolimits ({\mathbb X} , {\mathbb Y} )\} $.
Also let ${\mathcal I} _{{\mathbb X} }=\{ S\subset X: \neg \exists A\in {\mathbb P} ({\mathbb X} )\; A\subset S \}$.
We will use the following statements.
\begin{fac}[\cite{Ktow}]  \rm \label{T4057}
For each relational structure ${\mathbb X}$ we have:
$|\mathop{\rm sq}\nolimits  \langle  {\mathbb P} ({\mathbb X} ), \subset \rangle |\geq \aleph _0$ iff the poset $\langle  {\mathbb P} ({\mathbb X} ), \subset\rangle $ is atomless iff ${\mathbb P} ({\mathbb X} )$ contains two incompatible elements.
\end{fac}
\begin{fac}[\cite{Ktow}]  \rm \label{T4080}
A structure ${\mathbb X} $ is indivisible iff ${\mathcal I} _{{\mathbb X} }$ is an ideal in $P(X)$. Then

(a) $\mathop{\rm sm}\nolimits \langle  {\mathbb P} ({\mathbb X} ), \subset \rangle = \langle  {\mathbb P} ({\mathbb X} ), \subset _{{\mathcal I} _{\mathbb X}} \rangle $, where
$A \subset _{{\mathcal I} _{\mathbb X}} B \Leftrightarrow A\setminus B \in {\mathcal I} _{\mathbb X} $;

(b) $\mathop{\rm sq}\nolimits  \langle  {\mathbb P} ({\mathbb X} ), \subset \rangle $ is isomorphic to a dense subset of  $ \langle  (P (X )/\!\! =_{{\mathcal I} _{\mathbb X}})^+, \leq _{{\mathcal I} _{\mathbb X}} \rangle $. Hence
the poset $\langle  {\mathbb P} ({\mathbb X} ) , \subset \rangle $ is forcing equivalent to $(P(X)/{\mathcal I} _{{\mathbb X} })^+$.

(c) If ${\mathbb X}$ is countable, then $\langle  {\mathbb P} ({\mathbb X} ) , \subset \rangle $ is an atomless partial order of size ${\mathfrak c}$.
\end{fac}
\begin{fac}[\cite{Ktow}]  \rm\label{T4015}
Let ${\mathbb X} _i = \langle  X_i , \rho _i \rangle , i\in I$, and ${\mathbb Y} _j = \langle  Y_j , \sigma _j \rangle , j\in J$, be two families of disjoint connected
$L$-structures and ${\mathbb X}$ and ${\mathbb Y}$ their unions.
Then

(a) $F:{\mathbb X}  \hookrightarrow {\mathbb Y} $ iff $F=\bigcup _{i\in I} g_i$, where $f:I\rightarrow J$, $g_i : {\mathbb X} _i \hookrightarrow {\mathbb Y} _{f(i)}$, $i\in I$, and
\begin{equation}\label{EQ4087}
\forall \{i_1 ,i_2\} \in [I]^2 \;\; \forall x _{i_1} \in X_{i_1} \;\; \forall x _{i_2} \in X_{i_2} \;\;\neg \; g_{i_1}(x_{i_1})\; \sigma _{rs} \;g_{i_2}(x_{i_2});
\end{equation}

(b) $C\in {\mathbb P} ({\mathbb X} )$ iff $C= \bigcup _{i\in I} g_i [X_i]$, where $f:I\rightarrow I$, $g_i : {\mathbb X} _i \hookrightarrow {\mathbb X} _{f(i)}$, $i\in I$, and
\begin{equation}\label{EQ4009}
\forall \{i,j\} \in [I]^2 \;\; \forall x \in X_i \;\; \forall y \in X_j \;\;\neg \; g_i (x)\; \rho _{rs} \;g_j(y).
\end{equation}
\end{fac}
\begin{fac}[\cite{Ktow}]  \rm\label{T4026}
If ${\mathbb X}$ and ${\mathbb Y}$  are equimorphic structures, then the posets
$\langle  {\mathbb P} ({\mathbb X} ), \subset \rangle $ and $\langle  {\mathbb P} ({\mathbb Y} ), \subset \rangle $ are forcing equivalent.
\end{fac}
\begin{fac}[Pouzet \cite{Pouz}]  \rm  \label{T4094}
 If $p\leq |X|$ and ${\mathbb X}$ is $p$-monomorphic, then ${\mathbb X}$ is $r$-monomorphic for each $r\leq \min \{ p, |X|-p\}$. (See also \cite{Fra}, p.\ 259.)
\end{fac}
\section{Structures with maximally embeddable components}\label{S4006}
\begin{te}     \label{T4045}\rm
Let  ${\mathbb X} _i=\langle  X_i , \rho _{X_i}\rangle $, $i\in I$, be the components of a countable $L$-structure ${\mathbb X}=\langle  X, \rho \rangle $
and, for all $i,j\in I$, let

(i) ${\mathbb P} ({\mathbb X} _i , {\mathbb X} _j )=[{\mathbb X} _j ]^{|{\mathbb X} _i|}$ (the components of ${\mathbb X}$ are maximally embeddable),

(ii) $\forall A,B \in [{\mathbb X} _j ]^{|{\mathbb X} _i|}\;\; \exists a\in A \;\; \exists b\in B \;\; a\;\rho _{rs}\;b$.

\noindent
If $N=\{ |X_i |: i\in I \} $, $N_{\mathrm{fin}}=N\setminus \{ \omega \}$, $I_\kappa =\{ i\in I : |X_i |=\kappa \}$, for $\kappa \in N$,
$|I_\omega|=\mu$ and $Y=\bigcup _{i\in I\setminus I_\omega }X_i$, then we have

(a) $\mathop{\rm sq}\nolimits  \langle  {\mathbb P} ({\mathbb X} ), \subset \rangle $ is an $\omega _1$-closed atomless poset of size ${\mathfrak c}$. In addition, it is isomorphic
(and, hence, the poset $\langle  {\mathbb P} ({\mathbb X} ), \subset \rangle $ is forcing equivalent) to the poset
$$
                                            \begin{array}{llr}
                                              (P(\omega )/\mathop{\rm Fin}\nolimits )^+)^\mu     &  \mbox{ if } 1\leq \mu <\omega , \;  |N_{\mathrm{fin}}| <\omega  \mbox{ and } |Y|<\omega,&\mbox{ (a1)}\\
                                              ((P(\omega )/\mathop{\rm Fin}\nolimits )^+)^{\mu +1} &  \mbox{ if }0\leq \mu <\omega , \; |N_{\mathrm{fin}}| <\omega  \mbox{ and } |Y|=\omega,&\mbox{ (a2)}\\
                                              {\mathbb P} \times((P(\omega )/\mathop{\rm Fin}\nolimits )^+)^\mu  & \mbox{ if }0\leq \mu<\omega , \; |N_{\mathrm{fin}}| =\omega   ,&\mbox{ (a3)}\\                                               (P(\omega \times \omega)/(\mathop{\rm Fin}\nolimits  \times \mathop{\rm Fin}\nolimits  ))^+  &  \mbox{ if } \mu =\omega,&\mbox{ (a4)}
                                            \end{array}
$$
where ${\mathbb P}$ is an $\omega _1$-closed atomless poset, forcing
equivalent to $(P(\Delta )/{\mathcal E}{\mathcal D}_{\mathrm{fin}}  )^+ $.


(b) For some forcing related cardinal invariants of the poset $\langle  {\mathbb P} ({\mathbb X} ), \subset \rangle $ we have

\vspace{-3mm}
\begin{center}
{\scriptsize
\begin{tabular}{c|c|c|c}
If ${\mathbb X}$ satisfies         & $\langle  {\mathbb P} ({\mathbb X} ), \subset \rangle $ is     & $\mathop{\rm sq}\nolimits  \langle  {\mathbb P} ({\mathbb X} ), \subset \rangle $ is & ZFC $\vdash \mathop{\rm sq}\nolimits  \langle  {\mathbb P} ({\mathbb X} ), \subset \rangle $   \\
                          &  forcing equivalent to             &                                    & is ${\mathfrak h}$-distributive              \\[2mm] \hline
                          &                                    &                                    &                                              \\
$\mu < \omega \land |N_{\mathrm{fin}}|<\omega $  & $((P(\omega )/\mathop{\rm Fin}\nolimits  )^+)^n$, for some $n\in {\mathbb N}$            & ${\mathfrak t}$-closed          & yes iff $n=1$ \\[2mm]
$\mu < \omega \land |N_{\mathrm{fin}}|=\omega $  &  $(P(\Delta )/{\mathcal E}{\mathcal D}_{\mathrm{fin}}  )^+  \times  ((P(\omega )/\mathop{\rm Fin}\nolimits )^+)^\mu$ &  $\omega _1$-closed                & no            \\[2mm]
$\mu=\omega$                                 & $(P(\omega \times \omega)/(\mathop{\rm Fin}\nolimits  \times \mathop{\rm Fin}\nolimits  ))^+$             & $\omega _1$ but not $\omega _2$-closed & no
\end{tabular}
}
\end{center}
where $n=1$ iff $N\in [{\mathbb N} ]^{<\omega} \lor (|Y|<\omega \land \mu=1)$.

(c) ${\mathbb X}$ is indivisible iff $N\in [{\mathbb N} ]^\omega$  or $N=\{ 1 \}$ or $|I|=1$ or $|I_\omega |=\omega$.
\end{te}
A proof of the theorem, given at the end of this section, is based on the following five claims.
\begin{cla}  \rm \label{T4060}
$C\in {\mathbb P}({\mathbb X} )$ iff there is an injection $f:I\rightarrow I$  and there are $C_i \in [X_{f(i)}]^{|X_i|}$, $i\in I$, such that $C=\bigcup _{i\in I}C_i$.
\end{cla}
\noindent{\bf Proof.}
($\Rightarrow$)
Let $C\in {\mathbb P}({\mathbb X} )$. By Fact \ref{T4015}(b) there are functions $f:I\rightarrow I$ and $g_i : {\mathbb X} _i \hookrightarrow {\mathbb X} _{f(i)}$, $i\in I$,
satisfying (\ref{EQ4009}) and such that $C= \bigcup _{i\in I} g_i [X_i]$. By (\ref{EQ4009}) and (ii),
$f$ is an injection. Since $g_i : {\mathbb X} _i \hookrightarrow {\mathbb X} _{f(i)}$ we have
$C_i=g_i [X_i]\in {\mathbb P} ({\mathbb X} _i , {\mathbb X} _{f(i)}) =[{\mathbb X} _{f(i)} ]^{|{\mathbb X} _i|}$.

($\Leftarrow$) Suppose that $f$ and $C_i$, $i\in I$, satisfy the assumptions. Since
$[{\mathbb X} _{f(i)} ]^{|{\mathbb X} _i |} = {\mathbb P} ({\mathbb X} _i , {\mathbb X} _{f(i)} )$ there are $g_i : {\mathbb X} _i \hookrightarrow {\mathbb X} _{f(i)}$, $i\in I$,
such that $C_i=g_i [X_i]$. Since $f$ is an injection, for different $i,j\in I$ the sets $g_i[X_i]$ and $g_j[X_j]$
are in different components of ${\mathbb X}$ and, hence, we have (\ref{EQ4009}). By Fact \ref{T4015}(b), $C\in {\mathbb P}({\mathbb X} )$.
\hfill $\Box$ \par \vspace*{2mm}
\noindent
We continue the proof considering the following  cases and subcases.

 1. $N\subset {\mathbb N}$, with subcases $N\in [{\mathbb N}]^\omega$ (Claim \ref{T4061}) and $N\in [{\mathbb N} ]^{<\omega}$ (Claim \ref{T4062});

 2. $N\not\subset {\mathbb N}$, with subcases $|I_\omega|<\omega$ (Claim \ref{T4064}) and $|I_\omega |=\omega$ (Claim \ref{T4065}).

\vspace{2mm}
\noindent
{\bf Case 1:} $N\subset {\mathbb N}$.
\begin{cla}[Case 1.1]\rm \label{T4061}
If $N \in [{\mathbb N}]^\omega$, then

(a) ${\mathbb X}$ is an indivisible structure;

(b) $\mathop{\rm sq}\nolimits  \langle  {\mathbb P} ({\mathbb X} ), \subset \rangle $ is an $\omega _1$-closed  atomless poset;

(c) The structures ${\mathbb X}_i$, $i\in I $, are either full relations or complete graphs or reflexive or irreflexive linear orderings;

(d) There are structures ${\mathbb X}_n$, $n\in {\mathbb N} \setminus N$, such that $|X_n|=n$ and that the extended family
$\{ {\mathbb X} _i : i\in I \}\cup \{ {\mathbb X} _n : n\in {\mathbb N} \setminus N \}$ satisfies (i) and (ii);

(e) The poset $\langle  {\mathbb P} ({\mathbb X} ), \subset \rangle $ is forcing equivalent to $(P(\Delta )/{\mathcal E}{\mathcal D}_{\mathrm{fin}}  )^+$.
\end{cla}
\noindent{\bf Proof.}
Clearly, $N \in [{\mathbb N}]^\omega$ implies that $|I|=\omega$. First we prove
\begin{equation}\label{EQ4085}
S\in {\mathcal I} _{\mathbb X} \Leftrightarrow \exists n\in \omega \; \forall i\in I \; |S\cap X_i|\leq n .
\end{equation}

($\Rightarrow$)
Here, for convenience, we assume that $I=\omega$.
Suppose that for each $n\in \omega$ there is $i\in I$ such that $|S\cap X_i|> n $.
Then $I^S_{>n}=\{i\in \omega : |S\cap X_i|> n\}$, $n\in \omega$, are infinite sets.
By recursion we define sequences $\langle  i_k : k\in \omega \rangle $ in $\omega$ and $\langle  C_k : k\in \omega \rangle $ in $P(X)$ such that for each $k,l\in \omega$

(i) $k<l \Rightarrow i_k < i_l$,

(ii) $C_k \in [S\cap X_{i_k}]^{|X_k|}$.

\noindent
Suppose that the sequences  $i_0, \dots , i_k$ and $C_0 , \dots , C_k$ satisfy (i) and (ii).
Since $|I^S_{>|X_{k+1}|}|=\omega$ there is $i_{k+1}=\min \{i>i_k : |S\cap X_i|> |X_{k+1}|\}$ so
$|S\cap X_{i_{k+1}}|>|X_{k+1}|$, we choose $C_{k+1}\in [S\cap X_{i_{k+1}}]^{|X_{k+1}|}$ and the recursion works.

By (i) the function $f:I\rightarrow I$ defined by $f(k)=i_k$ is an injection.
By (ii) we have $C_k \in [X_{f(k)}]^{|X_k |}$ and, by Claim \ref{T4060}
$C=\bigcup _{k\in \omega }C_k \in {\mathbb P} ({\mathbb X} )$. Since $C\subset S$ we have $S\not\in {\mathcal I} _{\mathbb X} $.

($\Leftarrow$) Suppose that $C\in {\mathbb P} ({\mathbb X} )$, where $C\subset S$.
By Claim \ref{T4060} there are an injection $f:I\rightarrow I$  and $C_i \in [X_{f(i)}]^{|X_i|}$, $i\in I$, such that $C=\bigcup _{i\in I}C_i$.
For $n\in \omega$ there is $i_0 \in I$ such that $|X_{i_0}|>n$ and, hence, $C_{i_0}\in [X_{f(i_0)}]^{|X_{i_0}|}$, which implies $|X_{f(i_0)}\cap S|\geq |C_{i_0}|>n$.
(\ref{EQ4085}) is proved.

(a)
Suppose that $X=C\cup D$ is a partition, where $C,D\in {\mathcal I} _{\mathbb X}$. Then, by (\ref{EQ4085}), there are $m,n\in \omega$ such that
$|C\cap X_i |\leq m$ and $|D\cap X_i |\leq n$, for each $i\in I$. Hence for each $i\in I$ we have $|X_i| = |(X_i \cap C )\cup (X_i \cap D)|\leq m+n$, which is impossible since, by the assumption, $N\in [{\mathbb N} ]^\omega$.

(b)
By Facts \ref{T4042}(b) and (c)  it is sufficient to show that $\mathop{\rm sm}\nolimits \langle  {\mathbb P} ({\mathbb X} ), \subset \rangle $ is an $\omega _1$-closed and atomless pre-order.
Let $\mathop{\rm sm}\nolimits \langle  {\mathbb P} ({\mathbb X} ), \subset \rangle = \langle  {\mathbb P} ({\mathbb X} ), \leq \rangle $.
By Fact \ref{T4080} and (a) for each $A,B\in {\mathbb P} ({\mathbb X} )$ we have $A\leq B $ iff $A\setminus B \in {\mathcal I} _{\mathbb X}$ and, by (\ref{EQ4085}),
\begin{equation}\label{EQ4086}
A\leq B \Leftrightarrow \exists n\in {\mathbb N} \;\; \forall i\in I \;\; |A\setminus B \cap X_i|\leq n .
\end{equation}
Let $A_n \in {\mathbb P} ({\mathbb X} )$, $n\in \omega$, and $A_{n+1}\leq A_n$, for all $n\in \omega$.
We will find $A \in {\mathbb P} ({\mathbb X} )$ such that $A\leq A_n$, for all $n\in \omega$, that is, by (\ref{EQ4086}),
\begin{equation}\label{EQ4013}
\forall n\in \omega \;\;\exists m \in {\mathbb N} \;\; \forall i\in I \;\; |A\setminus A_n \cap X_i|\leq m .
\end{equation}
By recursion we define a sequence $\langle  i_r :r\in \omega \rangle $ in $I$ such that for each $r,s\in \omega$

(i) $r\neq s \Rightarrow i_r \neq i_s$,

(ii) $|A_0 \cap A_1 \cap \dots \cap A_r \cap X_{i_r}|>r$.

\noindent
First we choose $i_0$ such that $|A_0 \cap X_{i_0}|>0$.
Let the sequence $i_0, \dots , i_r$ satisfy (i) and (ii). For each $k\leq r $ we have $A_{k+1}\leq A_k$ and, by (\ref{EQ4086}), there is $m_k \in \omega$ such that
$\forall i\in I \;\; |A_{k+1}\setminus A_k \cap X_i|\leq m_k$. Thus
\begin{equation}\label{EQ4015}\textstyle
\forall i\in I \;\; \forall k\leq r \;\; |A_{k+1}\setminus A_k \cap X_i|\leq m_k .
\end{equation}
Since $A_{r+1}\in {\mathbb P} ({\mathbb X} )$ and $N\in [{\mathbb N} ]^\omega$,  by Claim \ref{T4060} the set
\begin{equation}\label{EQ4016}\textstyle
J_{r+1}=\{ i\in I : |A_{r+1}\cap X_i | > (\sum _{k\leq r}m_k )+r+1 \}
\end{equation}
is infinite and we choose
\begin{equation}\label{EQ4017}\textstyle
i_{r+1}\in J_{r+1}\setminus \{ i_0 , \dots i_r \} .
\end{equation}
Then (i) is true. Clearly,
$A_{r+1} \subset (\bigcap _{k=0}^{r+1} A_k )\cup \bigcup _{k=0}^{r}(A_{k+1}\setminus A_k )$ and, hence,
$A_{r+1}\cap X_{i_{r+1}} \subset (\bigcap _{k=0}^{r+1} A_k \cap X_{i_{r+1}})\cup \bigcup _{k=0}^{r}(A_{k+1}\setminus A_k \cap X_{i_{r+1}})$.
So, by (\ref{EQ4015})-(\ref{EQ4017})
$(\sum _{k\leq r}m_k )+r+1
< |A_{r+1}\cap X_{i_{r+1}}|
\leq |\bigcap _{k=0}^{r+1} A_k \cap X_{i_{r+1}}|+ \sum _{k\leq r}m_k $, which implies
$ |A_0 \cap \dots \cap A_r \cap A_{r+1}\cap X_{i_{r+1}}| >r+1$ and (ii) is true. The recursion works.

Let $S=\bigcup _{r\in \omega } (A_0 \cap A_1 \cap \dots \cap A_r \cap X_{i_r})$. By (i), (ii) and (\ref{EQ4085}) we have
$S\not\in{\mathcal I}_{\mathbb X}$ and, hence, there is $A\in {\mathbb P} ({\mathbb X} )$ such that $A\subset S$. We prove (\ref{EQ4013}). For $n\in \omega$ we have
$A\setminus A_n \subset S\setminus A_n \subset \bigcup _{r<n } (A_0 \cap A_1 \cap \dots \cap A_r \cap X_{i_r})\subset \bigcup _{r<n } X_{i_r}$,
thus $|A\setminus A_n |=m$, for some $m\in \omega$ and, hence, $|A\setminus A_n \cap X_i|\leq m$, for each $i\in I$.

So $\mathop{\rm sq}\nolimits  \langle  {\mathbb P} ({\mathbb X} ), \subset \rangle $ is $\omega _1$-closed.
By (a) and Facts \ref{T4080}(c) and \ref{T4042}(b) it is atomless.

(c)
Since $N\in [{\mathbb N} ]^\omega$, there are $i_0,i_1 \in I$ such that $|X_{i_0}|\geq 3$ and $|X_{i_1}|\geq |X_{i_0}|+3$.
By (i) we have ${\mathbb P} ({\mathbb X} _{i_0} , {\mathbb X} _{i_1} )=[{\mathbb X} _{i_1} ]^{|{\mathbb X} _{i_0}|}$ and, hence, the structure ${\mathbb X} _{i_1}$ is
$|X_{i_0}|$-monomorphic. Since $|X_{i_1}|- |X_{i_0}|\geq 3$ we have $\min \{ |X_{i_0}|,|X_{i_1}|- |X_{i_0}| \}\geq 3 $ and, by
Fact \ref{T4094},
\begin{equation}\label{EQ4070}
\forall r\leq 3\;\;({\mathbb X}_{i_1} \mbox{ is }r\mbox{-monomorphic}).
\end{equation}
Let $\{y_1 , y_2, y_3\}\in [X_{i_1}]^3$ and, for $r\in \{1,2,3\}$, let ${\mathbb Y}_r=\langle  Y_r , \tau  _r\rangle $, where $Y_r=\{ y_k :k\leq r\}$ and $\tau  _r =(\rho _{i_1})_{Y_r}$.
We prove
\begin{equation}\label{EQ4071}
\forall i\in I \;\; \forall r\leq \min \{ 3, |X_i |\} \;\;\forall A\in [X_i]^r \;\; \langle  A , (\rho _i )_A \rangle \cong {\mathbb Y} _r.
\end{equation}
If $|X_i|\geq |X_{i_1}|$, let $A\subset B\in [X_i]^{|X_{i_1}|}$. By (i) there exists an isomorphism
$f: \langle  B , (\rho _i)_B \rangle \rightarrow {\mathbb X} _{i_1}$ and, by (\ref{EQ4070}) we have
$\langle  A , (\rho _i )_A \rangle \cong \langle  f[A] , (\rho _{i_1} )_{f[A]} \rangle \cong {\mathbb Y} _r$.

If $|X_i|< |X_{i_1}|$ then, by (i), there exists an isomorphism
$f: {\mathbb X} _i \rightarrow {\mathbb X} _{i_1}$ and by (\ref{EQ4070}) we have
$\langle  A , (\rho _i )_A \rangle \cong \langle  f[A] , (\rho _{i_1} )_{f[A]} \rangle \cong {\mathbb Y} _r$.
Thus (\ref{EQ4071}) is true.

Clearly we have $\tau  _1=\emptyset$ or $\tau  _1=\{ \langle  y_1,y_1 \rangle \}$.

First, suppose that $\tau  _1=\emptyset$.
Then by (\ref{EQ4071}), for each $i\in I$ we have
\begin{equation}\label{EQ4072}
\forall x\in X_i \;\;\neg x\; \rho _i \; x,
\end{equation}
that is, all relations $\rho _i$, $i\in I$, are irreflexive. Suppose that $\tau  _2 \cap \{\langle  y_1, y_2 \rangle , \langle  y_2 , y_1 \rangle \}=\emptyset$.
Then by (\ref{EQ4071}) we would have $\rho _{i_1}=\emptyset$ and ${\mathbb X} _{i_1}$ would be a disconnected structure, which is not true. Thus
$\tau  _2 \cap \{\langle  y_1, y_2 \rangle , \langle  y_2 , y_1 \rangle \}\neq \emptyset$.

Thus, if $\langle  y_1, y_2 \rangle , \langle  y_2 , y_1 \rangle \in \tau  _2$, then by (\ref{EQ4071}), for each $i\in I$ we have
\begin{equation}\label{EQ4073}
\forall \{x,y \} \in [X_i ]^2 \;\; (x\; \rho _i \; y \;\land \; y\; \rho _i \; x)
\end{equation}
and, hence, ${\mathbb X}_i$ is a complete graph.

Otherwise, if $|\tau  _2 \cap \{\langle  y_1, y_2 \rangle , \langle  y_2 , y_1 \rangle \}|=1$ then, by (\ref{EQ4071}), for each $i\in I$ we have
\begin{equation}\label{EQ4074}
\forall \{x,y \} \in [X_i ]^2 \;\; (x\; \rho _i \; y \;\veebar \; y\; \rho _i \; x)
\end{equation}
and, hence, ${\mathbb X}_i$ is a tournament.
Thus ${\mathbb Y}_3$ is a tournament with three nodes and, hence,
${\mathbb Y} _3\cong C_3=\langle  \{1,2,3\} ,\{\langle  1,2 \rangle , \langle  2,3 \rangle , \langle  3,1 \rangle\}\rangle $ (the oriented circle graph)
or
${\mathbb Y} _3\cong L_3=\langle  \{1,2,3\} ,\{\langle  1,2 \rangle , \langle  2,3 \rangle , \langle  1,3 \rangle\}\rangle $ (the transitive triple, the strict linear order of size 3).
But ${\mathbb Y} _3\cong C_3$ would imply that ${\mathbb X} _{i_1}$ contains  a four element tournament having all substructures of size 3 isomorphic to
$C_3$, which is impossible. Thus ${\mathbb Y} _3\cong L_3$ which, together with (\ref{EQ4071}), (\ref{EQ4072}) and (\ref{EQ4074}) implies that all
relations $\rho _i$, $i\in I$ are transitive, so ${\mathbb X} _i$, $i\in I$, are strict linear orders.

If $\tau  _1=\{ \langle  y_1,y_1 \rangle\}$ then using the same arguments we show that the structures ${\mathbb X} _i$, $i\in I$, are either full relations or
reflexive linear orders.

(d)
follows from (c). Namely, if, for example, ${\mathbb X}_i$ are complete graphs, then ${\mathbb X} _n$ are complete graphs of size $n$.

(e) Let $N=\{ n_k : k\in {\mathbb N} \}$, where $n_1<n_2 < \dots $ and let ${\mathbb X}_n$, $n\in {\mathbb N} \setminus N$, be the structures from (d).
 W.l.o.g.\ suppose that
$I_{n_k} =\{ n_k \} \times \{ 1,2, \dots ,|I_{n_k}|\}$, if  $|I_{n_k}|\in {\mathbb N}$, and
$I_{n_k} =\{ n_k \} \times {\mathbb N} $, if  $|I_{n_k}|=\omega $.
Then $I\subset {\mathbb N} \times {\mathbb N}$ and $X=\bigcup _{k\in {\mathbb N} }\bigcup _{\langle  n_k ,r \rangle \in I_{n_k}} X_{\langle  n_k ,r \rangle }$.
For $l\in {\mathbb N}$, let ${\mathbb Y} _l =\langle  Y_l , \rho _l \rangle $ be defined by
$$
                           {\mathbb Y}_l=   \left\{ \begin{array}{ll}
                                                 {\mathbb X}_l                & \mbox{ if } l\in {\mathbb N} \setminus N,\\
                                                 {\mathbb X}_{\langle  n_k ,1 \rangle } & \mbox{ if } l=n_k, \mbox{ for a }k\in {\mathbb N} .
                                            \end{array}
                                    \right.
$$
and let ${\mathbb Y}=\langle  \bigcup _{l\in {\mathbb N}}Y_l ,\bigcup _{l\in {\mathbb N}}\rho _l \rangle $. We prove that ${\mathbb X} \hookrightarrow {\mathbb Y}$ and ${\mathbb Y} \hookrightarrow {\mathbb X}$.

${\mathbb Y} \hookrightarrow {\mathbb X}$. Let $f: {\mathbb N} \rightarrow I$, where $f(l)=\langle  n_l , 1\rangle $. Since $n_1<n_2 < \dots $ we have
$|Y_l|=l\leq n_l =|X_{\langle  n_l , 1\rangle }|=|X_{f(l)}|$ and, since the extended family of structures satisfies (i), there is
$g_l : {\mathbb Y} _l \hookrightarrow {\mathbb X} _{f(l)}$. Since $f$ is an injection, the sets $g_l [Y_l]$, $l\in {\mathbb N}$, are in different
components of ${\mathbb X}$ and, hence, condition (\ref{EQ4087})  is satisfied. Thus, by Fact \ref{T4015}(a), $F=\bigcup _{l\in {\mathbb N}}g_l :{\mathbb Y} \hookrightarrow {\mathbb X}$.

${\mathbb X} \hookrightarrow {\mathbb Y}$. Let ${\mathbb N} =\bigcup _{k\in {\mathbb N}}J_k$ be a partition, where $|J_k|=\omega$, for each $k\in {\mathbb N}$, and let
$Z_k=\bigcup _{\langle  n_k , r \rangle \in I_{n_k}}X_{\langle  n_k , r \rangle }$ and $T_k =\bigcup _{l\in J_k}Y_l$, for $k\in {\mathbb N}$.
Now $|I_{n_k}|\leq \omega =|J_k|$ and for $l\geq n_k$ we have $|X_{\langle  n_k , r \rangle }|=n_k\leq l =|Y_l|$. Hence
there is an injection $f_k : I_{n_k}\rightarrow J_k \setminus n_k$ and, since the extended family satisfies (i),
there are embeddings $g_{\langle  n_k , r \rangle }: {\mathbb X} _{\langle  n_k , r \rangle }\hookrightarrow {\mathbb Y} _{f(\langle  n_k , r \rangle )}$, for $\langle  n_k , r \rangle \in I_{n_k}$.
Thus, $f=\bigcup _{k\in {\mathbb N}} f_k : I \rightarrow {\mathbb N}$ and condition (\ref{EQ4087})  is satisfied so, by Fact \ref{T4015},
$F= \bigcup _{k\in {\mathbb N} }\bigcup _{\langle  n_k ,r \rangle \in I_{n_k}} g_{\langle  n_k ,r \rangle }$ embeds
${\mathbb X}=\bigcup _{k\in {\mathbb N} }\bigcup _{\langle  n_k ,r \rangle \in I_{n_k}} {\mathbb X}_{\langle  n_k ,r \rangle }$ into
${\mathbb Y}=\bigcup _{k\in {\mathbb N}}\bigcup _{l\in J_k}{\mathbb Y}_l$.

Now, by Fact \ref{T4026}, the posets $\langle  {\mathbb P} ({\mathbb X} ), \subset \rangle $ and $\langle  {\mathbb P} ({\mathbb Y} ), \subset \rangle $ are forcing equivalent.
W.l.o.g.\ suppose that $Y_l = \{ l \} \times \{ 1,2, \dots ,l \} \subset {\mathbb N} \times {\mathbb N}$. Then
$Y=\Delta = \{\langle  l,m \rangle \in {\mathbb N} \times {\mathbb N} : m\leq l \}$ and, by (\ref{EQ4085}), $S\in {\mathcal I} _{\mathbb Y}$ iff $\exists n\in {\mathbb N} \; \forall l\in {\mathbb N} \; |S\cap Y_l|\leq n$
iff $S\in {\mathcal E}{\mathcal D}_{\mathrm{fin}} $. Thus ${\mathcal I} _{\mathbb Y} ={\mathcal E}{\mathcal D}_{\mathrm{fin}} $ and, by Claim \ref{T4061}(a) and Fact \ref{T4080}(b),
$\langle  {\mathbb P} ({\mathbb Y} ), \subset \rangle $ is forcing equivalent to $(P(Y)/{\mathcal I} _{\mathbb Y} )^+$, that is to $(P(\Delta )/{\mathcal E}{\mathcal D}_{\mathrm{fin}}  )^+$.
\hfill $\Box$

\begin{cla}[Case 1.2]\rm \label{T4062}
If $N\in [{\mathbb N} ]^{<\omega }$, then we have

(a) $\mathop{\rm sq}\nolimits  \langle  {\mathbb P} ({\mathbb X} ), \subset  \rangle \cong (P(\omega )/\mathop{\rm Fin}\nolimits  )^+$;

(b) ${\mathbb X}$ is an indivisible structure iff $m=1$, where $m=\max N$.
\end{cla}
\noindent{\bf Proof.}
(a) Case A: $|I_m|=\omega$. For $S\subset X$ let $I_m^S=\{ i\in I_m : X_i \subset S\}$. First we prove
\begin{equation}\label{EQ4097}
S\in {\mathcal I} _{\mathbb X} \Leftrightarrow |I_m^S|<\omega .
\end{equation}
Let $S\not\in {\mathcal I} _{\mathbb X}$ and $C\subset S$, where $C\in {\mathbb P} ({\mathbb X})$. By Claim \ref{T4060}
there are an injection $f:I\rightarrow I$  and $C_i \in [X_{f(i)}]^{|X_i|}$, $i\in I$, such that $C=\bigcup _{i\in I}C_i$.
For $i\in I_m$ we have $|X_i|=m$ and, since  $C_i \in [X_{f(i)}]^m$, we have $|X_{f(i)}|=m$ and $C_i=X_{f(i)}\subset S$.
Thus $f(i)\in I_m^S$, for each $i\in I_m$ which, since $f$ is one-to-one, implies $|I_m^S|=\omega$.

Suppose that $|I_m^S|=\omega$ and let $f:I\rightarrow I_m^S$ be a bijection. For $i\in I$ we have $X_{f(i)}\subset S$ and
$|X_i| \leq m =|X_{f(i )}|$ and we choose $C_i \in [X_{f(i)}]^{|X_i|}$. Now  $C=\bigcup _{i\in I}C_i \subset S$ and, by Claim \ref{T4060},
$C\in {\mathbb P} ({\mathbb X})$. Thus $S\not\in {\mathcal I} _{\mathbb X}$ and (\ref{EQ4097}) is proved.

W.l.o.g.\ we assume that $I_m =\omega$.
By (\ref{EQ4097}), for $A\in {\mathbb P} ({\mathbb X} )$ we have $I^A_m \in [\omega ]^\omega$ and
we show that the posets $\langle  {\mathbb P} ({\mathbb X} ),\subset \rangle $ and $\langle  [\omega ]^\omega , \subset \rangle $
and the mapping $f:{\mathbb P} ({\mathbb X} )\rightarrow [\omega ]^\omega$
defined by $f(A)=I^A_m$ satisfy the assumptions of Fact \ref{T4090}.
Clearly, $A\subset B$ implies $I^A_m \subset I^B_m$ and (i) is true.
If $A$ and $B$ are incompatible elements of ${\mathbb P} ({\mathbb X})$, that is $A\cap B \in {\mathcal I} _{\mathbb X}$, then, by (\ref{EQ4097}), we have
$|I^{A\cap B}_m|<\omega$ and, since $I^{A}_m \cap I^{B}_m = I^{A\cap B}_m$, $f(A)$ and $f(B)$ are
incompatible in the poset $\langle  [\omega ]^\omega , \subset \rangle $. Thus (ii) is true as well.

We prove that $f$ is a surjection. Let $S\in [\omega ]^\omega$ and let $g:\omega \rightarrow S$ be a bijection.
Then $h=\mathop{\mbox{id}}\nolimits_{I\setminus \omega }\cup \;g :I\rightarrow I$ is an injection.
For $i\in \omega$ we have $h(i)=g(i)\in S$ and
we define $C_i=X_{g(i)}\in [X_{g(i)}]^{|X_i|}$. For $i\in I\setminus \omega$ let $C_i=X_i$.
Then, by Claim \ref{T4060}, $C=\bigcup _{i\in I}C_i=\bigcup _{i\in I\setminus \omega}X_i \cup \bigcup _{i\in \omega} X_{g(i)} \in {\mathbb P} ({\mathbb X} )$.
Now we have $f(C)= I^C_m = \{ g(i):i\in \omega \}=S$.

By Fact \ref{T4090}, $\mathop{\rm sq}\nolimits  \langle  {\mathbb P} ({\mathbb X} ), \subset  \rangle \cong \mathop{\rm sq}\nolimits  \langle  [\omega ]^\omega , \subset \rangle =(P(\omega )/\mathop{\rm Fin}\nolimits  )^+$.

Case B: $|I_m|<\omega$. Since $|X|=\omega$ the set $I=\bigcup _{n\in N}I_n$ is infinite and, hence,  there is $m_0=\max \{ n\in N : |I_n| =\omega \}$.
Clearly we have
\begin{equation}\label{EQ4018}
|I_{m_0}| =\omega \;\; \mbox{ and }\;\;\forall n\in N\setminus [0 , m_0] \;\;|I_n|<\omega
\end{equation}
and $X=Y\cup Z$, where $Y=\bigcup _{n\in N \cap[0, m_0]}\bigcup _{i\in I_n} X_i$ and $Z=\bigcup _{n\in N \setminus [0, m_0]}\bigcup _{i\in I_n} X_i$.
If $A\in {\mathbb P} ({\mathbb X})$, then for each $n\in N \setminus [0, m_0]$ the copy $A$ has exactly $|I_n|$-many components of size $n$ and, by (\ref{EQ4018}) and Claim \ref{T4060},
$Z\subset A$. So, it is easy to see that ${\mathbb P} ({\mathbb X} )= \{ C\cup Z: C\in {\mathbb P} ({\mathbb Y} )\}$ and, hence, the mapping
$F: {\mathbb P} ({\mathbb Y} ) \rightarrow {\mathbb P} ({\mathbb X} )$ given by $F(C)=C\cup Z$ is well defined and onto.
If $F(C_1)= F(C_2)$ then $(C_1 \cup Z) \cap Y = (C_2 \cup Z) \cap Y$, which implies $C_1=C_2$, thus $F$ is an injection. Clearly $C_1\subset C_2$ implies $F(C_1)\subset F(C_2)$ and,
if $F(C_1)\subset F(C_2)$, then $(C_1 \cup Z) \cap Y \subset (C_2 \cup Z) \cap Y$, which implies $C_1\subset C_2$. Thus
$\langle {\mathbb P} ({\mathbb X} ), \subset \rangle \cong _F  \langle {\mathbb P} ({\mathbb Y} ),\subset\rangle $ and,
by Fact \ref{T4042}(d), $\mathop{\rm sq}\nolimits  \langle {\mathbb P} ({\mathbb X} ), \subset \rangle \cong \mathop{\rm sq}\nolimits  \langle {\mathbb P} ({\mathbb Y} ),\subset\rangle $.
By (\ref{EQ4018}) the structure ${\mathbb Y}$ satisfies the assumption of Case A and, hence, $\mathop{\rm sq}\nolimits  \langle  {\mathbb P} ({\mathbb X} ), \subset \rangle \cong (P(\omega )/\mathop{\rm Fin}\nolimits  )^+$.

(b) If $m>1$, then there is a partition $X=A\cup B$ such that $A\cap X_i\neq \emptyset$ and $B\cap X_i\neq \emptyset$, for each $i\in I_m$.
Now, neither $A$ nor $B$ have a component of size $m$ and, hence, does not contain a copy of ${\mathbb X}$. Thus ${\mathbb X}$ is not indivisible.

If $m=1$, then $N=\{1\}$ and, since ${\mathbb P} ({\mathbb X} _i , {\mathbb X} _j )=[{\mathbb X} _j ]^{|{\mathbb X} _i|}$, the structures  ${\mathbb X}_i =\langle  \{ x_i \} , \rho _{\{ x_i \}} \rangle $, $i\in I$, are isomorphic and, hence, either $\rho _{\{ x_i \}} =\emptyset$, for all $i\in I$, which implies $\rho =\emptyset$
or $\rho _{\{ x_i \}} =\{ \langle  x_i,x_i \rangle\}$, for all $i\in I$, which implies $\rho =\Delta _X$. Thus, since $|I|=\omega$,
either ${\mathbb X} \cong \langle  \omega , \emptyset \rangle $ or ${\mathbb X} \cong \langle  \omega , \Delta _\omega \rangle $ and ${\mathbb P} ({\mathbb X} )= [X]^\omega $ in both cases,
which implies that ${\mathbb X}$ is an indivisible structure.
\hfill $\Box$ \par \vspace*{2mm}
\noindent
{\bf Case 2:} $N\not\subset {\mathbb N}$. Then $\mu >0$,
$X=(\bigcup _{i\in I\setminus I_\omega }X_i) \;\dot{\cup }\; (\bigcup _{i\in I_\omega }X_i )= Y\;\dot{\cup }\;Z$ (maybe $Y=\emptyset$)
and ${\mathbb X}$ is the disjoint union of the structures ${\mathbb Y} =\langle  Y, \rho _Y \rangle $ and ${\mathbb Z} =\langle  Z, \rho _Z \rangle $.
\begin{cla}[Case 2.1]\rm \label{T4064}
If $\mu \in {\mathbb N}$, then

(a)
\vspace{-5mm}
\begin{equation}\label{EQ4083}
\mathop{\rm sq}\nolimits  \langle  {\mathbb P} ({\mathbb X} ), \subset \rangle \cong \left\{ \begin{array}{cl}
                                                 ((P(\omega )/\mathop{\rm Fin}\nolimits )^+)^{\mu }     &  \mbox{ if } |N_{\mathrm{fin}}|<\omega  \mbox{ and } |Y|<\omega,\\
                                                 ((P(\omega )/\mathop{\rm Fin}\nolimits )^+)^{\mu +1} &  \mbox{ if } |N_{\mathrm{fin}}|<\omega  \mbox{ and } |Y|=\omega,\\
                                                 {\mathbb P} \times ((P(\omega )/\mathop{\rm Fin}\nolimits )^+)^\mu     &  \mbox{ if } |N_{\mathrm{fin}}|= \omega,
                                            \end{array}
                                    \right.
\end{equation}
where ${\mathbb P}$ is an $\omega _1$-closed atomless poset;

(b) If $|N_{{\mathrm{fin}}}|= \omega$, then  $\langle  {\mathbb P} ({\mathbb X} ), \subset \rangle $ and $(P(\Delta )/{\mathcal E}{\mathcal D}_{\mathrm{fin}}  )^+ \times (P(\omega )/\mathop{\rm Fin}\nolimits )^+)^\mu$
are forcing equivalent posets;

(c) ${\mathbb X}$ is indivisible iff $|I|=1$, that is $Y=\emptyset$ and $\mu=1$.
\end{cla}
\noindent{\bf Proof.}
(a)
For $i\in I_\omega$, let $A_i , B_i \in [X_i ]^\omega$ be disjoint sets,
$A=\bigcup _{i\in I\setminus I_\omega}X_i \cup \bigcup _{i\in I_\omega }A_i$ and
$B=\bigcup _{i\in I\setminus I_\omega}X_i \cup \bigcup _{i\in I_\omega }B_i$.
Then, by Claim \ref{T4060}, $A,B\in {\mathbb P} ({\mathbb X} )$ and, since $A\cap B$ does not contain infinite components,
we have $A\cap B \in {\mathcal I} _{\mathbb X}$. By Facts \ref{T4057} and \ref{T4042}(b), the posets $\langle  {\mathbb P} ({\mathbb X} ) , \subset \rangle $ and
$\mathop{\rm sq}\nolimits  \langle  {\mathbb P} ({\mathbb X} ) , \subset \rangle $ are atomless.

Concerning the closure properties of $\mathop{\rm sq}\nolimits  \langle  {\mathbb P} ({\mathbb X} ) , \subset \rangle $, first we prove the equality
\begin{equation}\label{EQ4099}
{\mathbb P} ({\mathbb X} )=\{ A\cup B : A\in {\mathbb P} ({\mathbb Y} ) \land B\in {\mathbb P} ({\mathbb Z} )\} .
\end{equation}
If $C\in {\mathbb P} ({\mathbb X} )$, then, by Claim \ref{T4060}, there is an injection $f:I\rightarrow I$  and there are $C_i \in [X_{f(i)}]^{|X_i|}$, $i\in I$, such that $C=\bigcup _{i\in I}C_i$. For $i\in I _\omega$ we have $C_i \in [X_{f(i)}]^\omega$ and, hence, $f(i)\in I_\omega$. Thus $f[I_\omega]\subset I_\omega$ and, since $f$ is one-to-one
and $I_\omega$ is finite, $f[I_\omega]= I_\omega$ and $f[I\setminus I_\omega] \subset I\setminus I_\omega$. Now we have $C=A\;\dot{\cup }\; B$, where
$A=\bigcup _{i\in I\setminus I_\omega }C_i \subset Y$ and $B= \bigcup _{i\in I_\omega }C_i \subset Z$. Clearly the structures ${\mathbb Y}$ and ${\mathbb Z}$ satisfy the
assumptions of Theorem \ref{T4045} and, since the restrictions
$f\upharpoonright I\setminus I_\omega :I\setminus I_\omega\rightarrow I\setminus I_\omega$ and
$f \upharpoonright I_\omega : I_\omega \rightarrow I_\omega$ are injections, by Claim \ref{T4060} we have $A\in {\mathbb P} ({\mathbb Y} )$ and $B\in {\mathbb P} ({\mathbb Z} )$.

Let $A\in {\mathbb P} ({\mathbb Y} )$ and $B\in {\mathbb P} ({\mathbb Z} )$. Since the structures ${\mathbb Y}$ and ${\mathbb Z}$ satisfy the
assumptions of Theorem \ref{T4045}, by Claim \ref{T4060} there are injections $g:I\setminus I_\omega\rightarrow I\setminus I_\omega $
and $h : I_\omega \rightarrow I_\omega$ and there are $C_i \in [X_{g(i)}]^{|X_i|}$, $i\in I\setminus I_\omega$, and  $C_i \in [X_{h(i)}]^{|X_i|}$,
$i\in I_\omega$, such that $A=\bigcup _{i\in I\setminus I_\omega}C_i$ and $B=\bigcup _{i\in I_\omega }C_i$. Now
$f=g\cup h :I\rightarrow I$ is an injection, $C_i \in [X_{f(i)}]^{|X_i|}$, for all $i\in I$, and, by Claim \ref{T4060},
$A\cup B=\bigcup _{i\in I}C_i \in {\mathbb P}({\mathbb X} )$. Thus (\ref{EQ4099}) is true.

Now we prove that
\begin{equation}\label{EQ4080}
\mathop{\rm sq}\nolimits  \langle  {\mathbb P} ({\mathbb X} ) , \subset \rangle \cong \mathop{\rm sq}\nolimits  \langle  {\mathbb P} ({\mathbb Y} ), \subset \rangle \times  \mathop{\rm sq}\nolimits  \langle  {\mathbb P} ({\mathbb Z} ),\subset \rangle .
\end{equation}
By (\ref{EQ4099}), the function $F: {\mathbb P} ({\mathbb Y} ) \times {\mathbb P} ({\mathbb Z} ) \rightarrow {\mathbb P} ({\mathbb X})$ given by $F(\langle  A,B \rangle )=A\cup B$
is well defined and onto and, clearly, it is a monotone injection. If $F(\langle  A,B \rangle ) \subset F(\langle  A',B' \rangle )$, then
$(A\cup B)\cap Y \subset  (A'\cup B')\cap Y $, that is $A\subset A'$ and, similarly, $B\subset B'$,
thus $\langle  A,B \rangle \leq \langle  A',B' \rangle $. So $F$ is an isomorphism and (\ref{EQ4080}) follows from (d) and (f) of Fact \ref{T4042}.

If $|N_{\mathrm{fin}}| <\omega $, then $|Y|<\omega$ implies $|{\mathbb P} ({\mathbb Y} )|=1$ and, hence,  $\mathop{\rm sq}\nolimits  \langle  {\mathbb P} ({\mathbb Y} ), \subset \rangle $ $ \cong 1$;
otherwise, if $|Y|=\omega$, then, by Claim \ref{T4062}, $\mathop{\rm sq}\nolimits  \langle  {\mathbb P} ({\mathbb Y} ), \subset \rangle \cong (P(\omega )/\mathop{\rm Fin}\nolimits  )^+$. So
\begin{equation}\label{EQ4081}
\mathop{\rm sq}\nolimits  \langle  {\mathbb P} ({\mathbb Y} ), \subset \rangle \cong \left\{ \begin{array}{cc}
                                                   1               &  \mbox{ if } |N_{\mathrm{fin}}| <\omega \mbox{ and } |Y|<\omega,\\
                                                   (P(\omega )/\mathop{\rm Fin}\nolimits )^+ &  \mbox{ if } |N_{\mathrm{fin}}| <\omega \mbox{ and } |Y|=\omega.
                                            \end{array}
                                    \right.
\end{equation}
By the assumption, for $i,j\in I_\omega$ we have ${\mathbb P} ({\mathbb X} _i, {\mathbb X} _j )=[X_j]^\omega$. Since $|I_\omega|<\omega$, by Claim \ref{T4060} we have
${\mathbb P} ({\mathbb Z} )=\{ \bigcup _{i\in I_\omega} C_i : \forall i\in I_\omega \; C_i \in [X_i]^\omega \}$
which implies $\langle  {\mathbb P} ({\mathbb Z}), \subset \rangle \cong \prod _{i\in I_\omega} \langle  [X_i]^\omega , \subset \rangle \cong\langle  [\omega ]^\omega , \subset\rangle ^\mu$.
Since $\mathop{\rm sq}\nolimits  \langle  [\omega ]^\omega , \subset\rangle =(P(\omega )/\mathop{\rm Fin}\nolimits  )^+$, by (d) and (f) of Fact \ref{T4042} we have
\begin{equation}\label{EQ4082}
\mathop{\rm sq}\nolimits  \langle  {\mathbb P} ({\mathbb Z}), \subset \rangle \cong ((P(\omega )/\mathop{\rm Fin}\nolimits  )^+) ^\mu.
\end{equation}
Now, for $|N_{\mathrm{fin}}| <\omega$  (\ref{EQ4083}) follows from (\ref{EQ4080}), (\ref{EQ4081}) and (\ref{EQ4082}).
If $|N_{\mathrm{fin}}| =\omega$, then, by Claim \ref{T4061}, ${\mathbb P} =\mathop{\rm sq}\nolimits  \langle  {\mathbb P} ({\mathbb Y} ), \subset \rangle $ is $\omega _1$-closed atomless
and (\ref{EQ4083}) follows from (\ref{EQ4080}) and (\ref{EQ4082}).

(b) By Claim \ref{T4061}(e) and Fact \ref{T4042}(a),
the posets $\langle  {\mathbb P} ({\mathbb Y} ) , \subset \rangle $, $\mathop{\rm sq}\nolimits \langle  {\mathbb P} ({\mathbb Y} ) , \subset \rangle $ and  $(P(\Delta )/{\mathcal E}{\mathcal D}_{\mathrm{fin}}  )^+$  are forcing equivalent.
By (\ref{EQ4080}) and (\ref{EQ4082}) we have
$\mathop{\rm sq}\nolimits \langle  {\mathbb P} ({\mathbb X} ) , \subset \rangle \cong \mathop{\rm sq}\nolimits \langle  {\mathbb P} ({\mathbb Y} ) , \subset \rangle  \times  (P(\omega )/\mathop{\rm Fin}\nolimits )^+)^\mu$.

(c) Let $Y=\emptyset$ and $\mu=1$. Then ${\mathbb P}({\mathbb X} )=[X]^\omega$ and, clearly, ${\mathbb X}$ is indivisible.

If $Y\neq\emptyset$, then, by (a), each $C\in {\mathbb P}({\mathbb X} )$ must intersect both $Y$ and $Z$ and the partition $X=Y\cup Z$ witnesses that
${\mathbb X}$ is not indivisible.

If $Y=\emptyset$ but $\mu>1$, by (a), each $C\in {\mathbb P}({\mathbb X} )$ must intersect all components of ${\mathbb X}$ and for
$i_0\in I_\omega =I$, the partition $X=X_{i_0} \cup \bigcup _{i\in I_\omega \setminus \{ i_0\}}X_i$ witnesses that
${\mathbb X}$ is not indivisible.
\hfill $\Box$
\begin{cla}[Case 2.2]\rm \label{T4065}
If $\mu=\omega$, then

(a) ${\mathbb X}$ is an indivisible structure;

(b) $\mathop{\rm sq}\nolimits  \langle  {\mathbb P} ({\mathbb X} ), \subset \rangle \cong (P(\omega \times \omega)/(\mathop{\rm Fin}\nolimits  \times \mathop{\rm Fin}\nolimits  ))^+$.
\end{cla}
\noindent{\bf Proof.}
(a) For $S\subset X$ let $I_\omega ^S =\{ i\in I_\omega : |S\cap X_i|=\omega \}$ and first we prove
\begin{equation}\label{EQ4098}
S\in {\mathcal I} _{\mathbb X} \Leftrightarrow |I_\omega^S |<\omega .
\end{equation}
Suppose that $|I_\omega^S |=\omega$. Let $f:I\rightarrow I_\omega^S$ be a bijection. Then, for $i\in I$ we have
$|S\cap X_{f(i)}|=\omega$ and we can choose $C_i \in [S\cap X_{f(i)}]^{X_i} \subset {\mathbb P} ({\mathbb X} _i , {\mathbb X} _{f(i)})$. By Claim
\ref{T4060} we have $C=\bigcup _{i\in I}C_i \in {\mathbb P} ({\mathbb X} )$ and, clearly, $C\subset S$. Thus  $S\not\in {\mathcal I} _{\mathbb X} $.

Let $S\not\in {\mathcal I} _{\mathbb X} $ and $C\in {\mathbb P} ({\mathbb X} )$, where $C \subset S$.
By Claim \ref{T4060} there are an injection $f:I\rightarrow I$
and $C_i \in [X_{f(i)}]^{|X_i|}$, $i\in I$, such that $C=\bigcup _{i\in I}C_i$.
For $i\in I_\omega$ we have $C_i \in [X_{f(i)}]^\omega$, which implies $|S\cap X_{f(i)}|=\omega$, that is
$f(i)\in I_\omega ^S$. Thus $f[I_\omega] \subset I_\omega ^S$ and, since $f$ is one-to-one and $|I_\omega|=\omega$, we have $|I_\omega^S |=\omega$ and (\ref{EQ4098}) is proved.

Suppose that ${\mathbb X}$ is divisible and $X=A\cup B$, where $A,B\in {\mathcal I} _{\mathbb X}$.
Then, by (\ref{EQ4098}), $|I_\omega^A \cup I_\omega^B|<\omega $ and there is $i\in I_\omega \setminus (I_\omega^A \cup I_\omega^B)$.
Now, $|A\cap X_i|, |B\cap X_i|<\omega$, which is impossible since $X_i =(A\cap X_i) \cup (B\cap X_i)$ is an infinite set.

(b)
W.l.o.g.\ we suppose that $I_\omega =\omega$ and $X_i = \{ i \} \times \omega$, for $i\in \omega$. Then
$X=Y\cup (\omega \times \omega)$, where $Y=\bigcup _{i\in I\setminus \omega}X_i$.
Clearly, for $S\subset \omega\times\omega$,
\begin{equation}\label{EQ4096}
S\in\mathop{\rm Fin}\nolimits \times\mathop{\rm Fin}\nolimits  \Leftrightarrow|I^S_\omega|<\omega.
\end{equation}
By (\ref{EQ4098}), for $A\in {\mathbb P} ({\mathbb X})$ the set $I_\omega ^A =I_\omega ^{A\cap (\omega \times \omega)}$
is infinite and by (\ref{EQ4096}) we have $A\cap (\omega \times \omega)\not\in  \mathop{\rm Fin}\nolimits  \times \mathop{\rm Fin}\nolimits $.
Hence the mapping
$$
f:\langle  {\mathbb P} ({\mathbb X} ), \subset \rangle \rightarrow \langle  (P(\omega \times \omega)/ \! =_{\mathop{\rm Fin}\nolimits  \times \mathop{\rm Fin}\nolimits  })^+ , \trianglelefteq _{\mathop{\rm Fin}\nolimits  \times \mathop{\rm Fin}\nolimits }\rangle $$
given by $f(A)=[A\cap (\omega \times \omega)]_{=_{\mathop{\rm Fin}\nolimits  \times \mathop{\rm Fin}\nolimits  }}$, for all $A\in {\mathbb P} ({\mathbb X} )$,
is well defined and we show that it satisfies the assumptions of Fact \ref{T4090}.
Let $A,B\in {\mathbb P} ({\mathbb X} )$.

(i) If $A\subset B$, then $(A\cap (\omega \times \omega))\setminus (B\cap (\omega \times \omega))=\emptyset \in \mathop{\rm Fin}\nolimits  \times \mathop{\rm Fin}\nolimits $
and $f(A)=[A\cap (\omega \times \omega)]_{=_{\mathop{\rm Fin}\nolimits  \times \mathop{\rm Fin}\nolimits  }} \trianglelefteq _{\mathop{\rm Fin}\nolimits  \times \mathop{\rm Fin}\nolimits } [B\cap (\omega \times \omega)]_{=_{\mathop{\rm Fin}\nolimits  \times \mathop{\rm Fin}\nolimits  }}=f(B)$.

(ii) If $A$ and $B$ are incompatible in  $\langle  {\mathbb P} ({\mathbb X} ), \subset \rangle $, then $A\cap B\in {\mathcal I} _{\mathbb X}$ and, by (\ref{EQ4098}), $|I^{A\cap B}_\omega|<\omega$, that is
$|I^{(A\cap (\omega \times \omega))\cap (B\cap (\omega \times \omega))}_\omega|<\omega$, which, by (\ref{EQ4096}) implies  $(A\cap (\omega \times \omega))\cap (B\cap (\omega \times \omega))\in  \mathop{\rm Fin}\nolimits  \times \mathop{\rm Fin}\nolimits $.
Hence $f(A)=[A\cap (\omega \times \omega)]_{=_{\mathop{\rm Fin}\nolimits  \times \mathop{\rm Fin}\nolimits  }}$ and $f(B)=[B\cap (\omega \times \omega)]_{=_{\mathop{\rm Fin}\nolimits  \times \mathop{\rm Fin}\nolimits  }}$ are incompatible
in $(P(\omega \times \omega)/ \! =_{\mathop{\rm Fin}\nolimits  \times \mathop{\rm Fin}\nolimits  })^+$.

(iii) We show that $f$ is a surjection.
It is easy to see that for $A,B\in {\mathbb P} ({\mathbb X} )$,
\begin{equation}\label{EQ4026}
I_\omega ^{A\setminus B}\cup I_\omega ^{B\setminus A} = I_\omega ^{A\Delta B}.
\end{equation}
Let $[S]_{=_{\mathop{\rm Fin}\nolimits  \times \mathop{\rm Fin}\nolimits  }} \in (P(\omega \times \omega)/ \! =_{\mathop{\rm Fin}\nolimits  \times \mathop{\rm Fin}\nolimits  })^+ $. Then, by (\ref{EQ4096}), we have  $|I^S_\omega|=\omega$.
 Let $g:\omega \rightarrow I_\omega ^S$ be a bijection.
Then $h=\mathop{\mbox{id}}\nolimits_{I\setminus \omega }\cup \;g :I\rightarrow I$ is an injection. For $i\in \omega$ we have $h(i)=g(i)\in I_\omega ^S$ and
we define $C_i=S\cap X_{g(i)}\in [X_{g(i)}]^{|X_i|}$. For $i\in I\setminus \omega$ let $C_i=X_i$. Then,
by Claim \ref{T4060},
$$\textstyle
C=\bigcup _{i\in I}C_i=\bigcup _{i\in I\setminus \omega}X_i \cup \bigcup _{i\in \omega} S\cap X_{g(i)} \in {\mathbb P} ({\mathbb X} ).
$$
Now $S\setminus C =\bigcup _{j\in \omega \setminus I^S_\omega}S\cap X_j$, which implies $I^{S\setminus C}_\omega =\emptyset$ and
$C\setminus S =\bigcup _{i\in I\setminus \omega}X_i \setminus S $, which implies $I^{C\setminus S}_\omega =\emptyset$.
So, by (\ref{EQ4026}), $I^{C\bigtriangleup S}_\omega =I^{(C\cap (\omega \times \omega))\bigtriangleup S}_\omega =\emptyset$
and, by (\ref{EQ4096}), $(C\cap (\omega \times \omega))\bigtriangleup S\in \mathop{\rm Fin}\nolimits \times \mathop{\rm Fin}\nolimits $,
so $f(C) = [C\cap (\omega \times \omega)]_{=_{\mathop{\rm Fin}\nolimits  \times \mathop{\rm Fin}\nolimits  }}=[S]_{=_{\mathop{\rm Fin}\nolimits  \times \mathop{\rm Fin}\nolimits  }}$.

By Fact \ref{T4090} and since $\langle  (P(\omega \times \omega)/ \! =_{\mathop{\rm Fin}\nolimits  \times \mathop{\rm Fin}\nolimits  })^+ , \trianglelefteq _{\mathop{\rm Fin}\nolimits  \times \mathop{\rm Fin}\nolimits }\rangle $ is a separative partial order we have
$\mathop{\rm sq}\nolimits  \langle  {\mathbb P} ({\mathbb X} ), \subset \rangle \cong \mathop{\rm sq}\nolimits  \langle  (P(\omega \times \omega)/ \! =_{\mathop{\rm Fin}\nolimits  \times \mathop{\rm Fin}\nolimits  })^+ , \trianglelefteq _{\mathop{\rm Fin}\nolimits  \times \mathop{\rm Fin}\nolimits }\rangle \cong \langle  (P(\omega \times \omega)/ \! =_{\mathop{\rm Fin}\nolimits  \times \mathop{\rm Fin}\nolimits  })^+ , \trianglelefteq _{\mathop{\rm Fin}\nolimits  \times \mathop{\rm Fin}\nolimits }\rangle $.
\hfill $\Box$ \par \vspace*{2mm}

\noindent
{\bf Proof of Theorem \ref{T4045}.}
(a)
(a4) is Claim \ref{T4065}(b).  For $\mu >0$, (a1)-(a3) are proved in Claim \ref{T4064}(a). For $\mu =0$, (a2) is proved in Claim \ref{T4062}(a)
and (a3) in Claim \ref{T4061}(b). By Facts \ref{T4043} and \ref{T4091},
$\mathop{\rm sq}\nolimits  \langle  {\mathbb P} ({\mathbb X} ), \subset \rangle $ is an $\omega _1$-closed atomless poset. It is of size ${\mathfrak c}$ since it contains
a reversed binary tree of height $\omega$  and the set of lower bounds of its branches is of cardinality ${\mathfrak c}$.
The forcing equivalent of ${\mathbb P}$ is given in Claim \ref{T4061}(e).

(b) follows from (a), Claim \ref{T4064}(b) and Fact \ref{T4091}.

(c)
The implication ``$\Leftarrow$" follows from  Claims \ref{T4061}(a), \ref{T4062}(b), \ref{T4064}(c)  and \ref{T4065}(a).
For a proof of ($\Rightarrow$) suppose that $N\not\in [{\mathbb N} ]^\omega$, $N\neq \{ 1 \}$, $|I|\neq 1$ and $|I_\omega |<\omega$.

If $N\subset {\mathbb N}$, then, since $N\not\in [{\mathbb N} ]^\omega$, we have $N=\{ n_0, \dots ,n_m \}$, where $n_0< \dots < n_m$ and, since $N\neq \{ 1 \}$, $n_m >1$.
Let $x_i \in X_i$, for $i\in I_{n_m}$, let
$A=\bigcup _{i\in I\setminus I_{n_m}} X_i \cup \bigcup _{i\in I_{n_m}} \{ x_i \}$ and
$B=\bigcup _{i\in I_{n_m}} X_i \setminus \{ x_i \}$. Then $X=A\cup B$ and neither $A$ nor $B$ contain a copy of ${\mathbb X}$, since all their components
are of size $<n_m$.

If $N\not\subset {\mathbb N}$, then $I_\omega \neq \emptyset $ and, since  $|I_\omega |<\omega$, we have
$0<|I_\omega |=m\in {\mathbb N}$. Since $|I|\neq 1$, by Claim \ref{T4064}(c) ${\mathbb X}$ is not indivisible.
\hfill $\Box$
\section{Examples}\label{S4008}
\begin{ex}   \rm \label{EX4020}
Equivalence relations on countable sets.
If ${\mathbb X}=\langle  X ,\rho \rangle $, where $\rho $ is an equivalence relation on a countable set $X$, then, clearly, the components $X_i$, $i\in I$, of ${\mathbb X}$ are the equivalence classes determined by $\rho $ and
for each $i\in I$ the restriction $\rho _{X_i}$ is the full relation on $X_i$, which implies that
conditions (i) and (ii) of Theorem \ref{T4045} are satisfied.
Thus the poset $\mathop{\rm sq}\nolimits  \langle  {\mathbb P} ({\mathbb X} ), \subset \rangle $ is $\omega _1$-closed and atomless and, hence, ${\mathbb X}$ belongs to the column $D$ of Diagram \ref{F4001}. Some examples of such structures are given in Diagram \ref{F4002}, where $\bigcup _m F_n$ denotes the disjoint union of $m$ full relations on a set of size $n$.
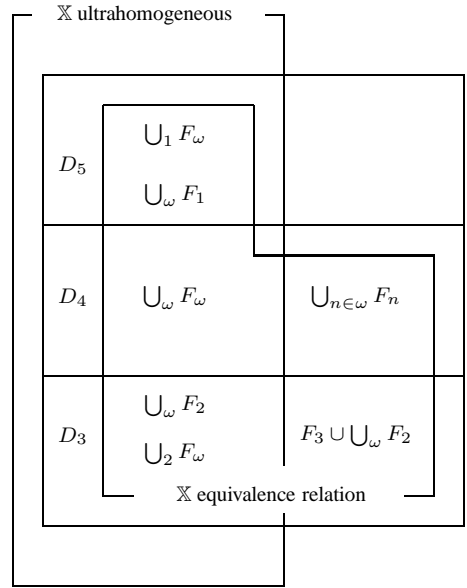
\begin{figure}[h]
\begin{center}
\unitlength  0.8mm 
\linethickness{0.4pt}
\ifx\plotpoint\undefined\newsavebox{\plotpoint}\fi 


\begin{picture}(85,102)(0,0)


\put(5,5){\line(1,0){45}}
\put(10,15){\line(1,0){70}}
\put(20,20){\line(1,0){5}}
\put(70,20){\line(1,0){5}}
\put(10,40){\line(1,0){70}}
\put(45,60){\line(1,0){30}}
\put(10,65){\line(1,0){70}}
\put(20,85){\line(1,0){25}}
\put(10,90){\line(1,0){70}}
\put(5,100){\line(1,0){3}}
\put(47,100){\line(1,0){3}}
\put(5,5){\line(0,1){95}}
\put(10,15){\line(0,1){75}}
\put(20,20){\line(0,1){65}}
\put(45,60){\line(0,1){25}}
\put(50,25){\line(0,1){75}}
\put(50,5){\line(0,1){12}}
\put(75,20){\line(0,1){40}}
\put(80,15){\line(0,1){75}}



\scriptsize

\put(27,100){\makebox(0,0)[cc]{${\mathbb X}$ ultrahomogeneous}}
\put(48,20){\makebox(0,0)[cc]{${\mathbb X} $ equivalence relation}}
\put(32,80){\makebox(0,0)[cc]{$\bigcup _ 1F_\omega$}}
\put(32,70){\makebox(0,0)[cc]{$\bigcup _\omega F_1$}}
\put(32,53){\makebox(0,0)[cc]{$\bigcup _\omega F_\omega$}}
\put(62,53){\makebox(0,0)[cc]{$\bigcup _{n\in \omega } F_n$}}
\put(32,35){\makebox(0,0)[cc]{$\bigcup _\omega  F_2$}}
\put(32,27){\makebox(0,0)[cc]{$\bigcup _2  F_\omega$}}
\put(62,30){\makebox(0,0)[cc]{$F_3 \cup \bigcup _\omega  F_2$}}
\put(15,30){\makebox(0,0)[cc]{$D_3$}}
\put(15,53){\makebox(0,0)[cc]{$D_4$}}
\put(15,75){\makebox(0,0)[cc]{$D_5$}}
\end{picture}

\end{center}

\vspace{-7mm}
\caption{Equivalence relations on countable sets}\label{F4002}
\end{figure}
We note that ${\mathbb X} $ is a ultrahomogeneous structure iff all equivalence classes are of the same size, so
the following countable equivalence relations are ultrahomogeneous and by Theorem \ref{T4045} have the given properties.

$\bigcup _{\omega }F_n$. It is indivisible iff $n=1$ (the diagonal) and the poset $\mathop{\rm sq}\nolimits  \langle  {\mathbb P} ({\mathbb X} ), \subset \rangle $ is isomorphic to
$(P(\omega )/\mathop{\rm Fin}\nolimits  )^+ $ which is a ${\mathfrak t}$-closed and ${\mathfrak h}$-distributive poset.

$\bigcup _{n}F_\omega$. It is indivisible iff $n=1$ (the full relation) and the poset $\mathop{\rm sq}\nolimits  \langle  {\mathbb P} ({\mathbb X} ), \subset \rangle $ is isomorphic to
$((P(\omega )/\mathop{\rm Fin}\nolimits  )^+)^n$  which is ${\mathfrak t}$-closed, but for $n>1$ not ${\mathfrak h}$-distributive poset in, for example, the Mathias model.

$\bigcup _{\omega }F_\omega $ (the $\omega$-homogeneous-universal equivalence relation). It is indivisible and
$\mathop{\rm sq}\nolimits  \langle  {\mathbb P} ({\mathbb X} ), \subset \rangle $ is isomorphic to $(P(\omega \times \omega)/(\mathop{\rm Fin}\nolimits  \times \mathop{\rm Fin}\nolimits  ))^+$, which is $\omega _1$-closed, but not $\omega _2$-closed and, hence,
consistently neither ${\mathfrak t}$-closed nor ${\mathfrak h}$-distributive.
\end{ex}
\begin{ex}   \rm \label{EX4021}
Disjoint unions of complete graphs.
The same picture as in Example \ref{EX4020} is obtained for
countable graphs ${\mathbb X} =\bigcup _{i\in I}{\mathbb X}_i$, where ${\mathbb X} _i =\langle  X_i ,\rho _i\rangle $, $i\in I$, are disjoint complete graphs
(that is $\rho _i=(X_i \times X_i) \setminus \Delta _{X_i}$) since, clearly, conditions (i) and (ii) of Theorem \ref{T4045} are satisfied. Also, by a well known
characterization of Lachlan and Woodrow \cite{Lach} all disconnected countable ultrahomogeneous graphs are of the form
$\bigcup _m K_n$ (the union of $m$-many complete graphs of size $n$), where $mn=\omega$ and $m>1$. So in Diagram \ref{F4002} we can replace $F_n$ with $K_n$.
\end{ex}
\begin{ex}   \rm \label{EX4022}
Disjoint unions of ordinals $\leq \omega$. A similar picture is obtained for
countable partial orders ${\mathbb X} =\bigcup _{i\in I}{\mathbb X}_i$, where ${\mathbb X} _i$'s are disjoint copies of ordinals $\alpha _i \leq \omega$.
(Clearly, linear orders satisfy (ii) of Theorem \ref{T4045} and ${\mathbb P} (\alpha , \beta )=[\beta ]^{|\alpha |}$, for each two ordinals $\alpha ,\beta \leq \omega$.)
So in Diagram \ref{F4002} we can replace $F_n$ with $L_n$, where $L_n\cong n\leq \omega$, but these partial orderings are not ultrahomogeneous.
\end{ex}
\begin{rem}   \rm \label{R4000}
All structures analyzed in Examples \ref{EX4020}, \ref{EX4021} and \ref{EX4022} are disconnected.
But, since ${\mathbb P} (\langle  X, \rho \rangle )={\mathbb P} (\langle  X, \rho ^c \rangle ) $,
taking their complements we obtain
connected structures with the same posets $\langle  {\mathbb P} ({\mathbb X} ), \subset\rangle $ and $\mathop{\rm sq}\nolimits \langle  {\mathbb P} ({\mathbb X} ), \subset\rangle $, having the properties
established in these examples. For example, the complement of $\bigcup _m F_n$ is the graph-theoretic complement of the graph $\bigcup _m K_n$.
\end{rem}
\begin{rem}   \rm \label{R4001}
The structures satisfying the assumptions of Theorem \ref{T4045}.
Let a countable structure ${\mathbb X} =\bigcup _{i\in I}{\mathbb X}_i$ satisfy conditions (i) and (ii).

First, (i) implies that all components of the same size are isomorphic.

Second, if $|X_i|=\omega$ for some $i\in I$, then, by (i),  ${\mathbb P} ({\mathbb X} _i )=[X _i ]^\omega$ and, by \cite{Ktow},
${\mathbb X} _i$ is isomorphic to one of the following structures:
1. The empty relation;
2. The complete graph;
3. The natural strict linear order on $\omega$;
4. Its inverse;
5. The diagonal relation;
6. The full relation;
7. The natural reflexive linear order on $\omega$;
8. Its inverse.
Thus, since ${\mathbb X} _i$ is a connected structure, it is isomorphic to the structure
2, 3, 4, 6, 7 or 8 and, by (i) again, this fact implies that
\begin{center}
$(\ast)$  All  ${\mathbb X} _i$'s are either full relations or complete graphs or linear orders.
\end{center}
By Claim \ref{T4061}(c), $(\ast)$ holds when ${\mathbb X}_i$'s are finite, but their sizes are unbounded.

But, if the size of the components of ${\mathbb X}$ is bounded by some $n\in {\mathbb N}$, there are structures
which do not satisfy $(\ast)$. For example, take a disjoint union of $\omega$ copies of the linear graph $L_n$ and
$\omega$ copies of the circle graph $C_{n+1}$.
\end{rem}

\end{document}